\documentclass[a4paper,oneside,10pt]{article}%
\usepackage[dvips]{color}
\usepackage[dvips]{graphics}
\usepackage[leqno]{amsmath}
\usepackage{amsmath}
\usepackage{amsthm}
\usepackage{amsfonts}
\usepackage{amssymb,latexsym}
\usepackage{amscd}
\usepackage{delarray}
\usepackage{enumerate}
\usepackage{graphicx}
\usepackage{epsfig}
\usepackage{amssymb}
\providecommand{\U}[1]{\protect\rule{.1in}{.1in}}

\def\B{{\mathcal{B}}}

\newtheorem{theorem}{Theorem}[section]

\newtheorem{definition}{Definition}[section]

\newtheorem{prop}{Proposition}[section]

\newtheorem{rque}{Remark}[section]
\begin{document}

\title{Asymptotic expansion of the solution of the steady Stokes equation with
variable viscosity in a two-dimensional tube structure}
\author{G.Cardone\\University of Sannio, Department of Engineering\\Piazza Roma, 21, 84100 Benevento, Italy\\email: giuseppe.cardone@unisannio.it
\and R.Fares, G.P.Panasenko\\Laboratory of Mathematics of the University of Saint Etienne (LaMUSE), EA 3989 \\23, rue P.Michelon, 42023 St. Etienne, France\\email: roula.fares@univ-st-etienne.fr; grigory.panasenko@univ-st-etienne.fr}
\maketitle

\begin{abstract}
The Stokes equation with the varying viscosity is considered in a thin tube
structure, i.e. in a connected union of thin rectangles with heights of order
$\varepsilon<<1 $ and with bases of order $1$ with smoothened boundary. An
asymptotic expansion of the solution is constructed: it contains some
Poiseuille type flows in the channels (rectangles) with some boundary layers
correctors in the neighborhoods of the bifurcations of the channels. The
estimates for the difference of the exact solution and its asymptotic
approximation are proved.

\medskip



\end{abstract}

\section{Introduction}

{The blood circulation, the transport of cells and substances in the human
body as well as some liquid-cooling systems and oil-recovery/oil-transport
processes in engineering, are modeled by the equations of fluid motion posed
in thin domains. The present paper studies the Stokes equation with varying
viscosity in a tube structure. In two-dimensional case a tube
structure\footnote{see Fig.\ref{f3}} (or pipe-wise structure) is some
connected union of thin rectangles with heights of order $\varepsilon$ and
with bases of order $1$ with a smoothened boundary $(\cite{CRAS}%
,\cite{pana2005})$; here $\varepsilon$ is a small positive parameter. An
asymptotic expansion for the case of a constant viscosity have been
constructed in \cite{BGP}. It is "compiled" of some Poiseuille flows inside
the rectangles glued by some boundary layer solutions, in the neighborhood of
the junctions. Here we consider a more general case when the viscosity is not
constant but depends on a longitudinal variable for each rectangle. This
situation models a blood flow in a vessel structure where the viscosity
depends on the concentration of some substances diluted in blood or some blood
cells. Indeed, the asymptotic analysis of the convection-diffusion equation
set in such domains $(\cite{CaPaSi},\cite{pana2005})$ shows that in the case
of the Neumann (impermeability) condition at the lateral boundary and small
Reynolds numbers, the concentration is asymptotically close to the
one-dimensional description, that is the convection-diffusion equation set on
the graph. The solution of the problem on the graph is the leading term of the
asymptotic expansion, and it evidently on depends on the longitudinal
variable. On the other hand, the viscosity often depends on the concentration
of the diluted substances or distributed cells, and so, it depends on the
longitudinal variable. Of course, the fluid motion equation is coupled with
the diffusion-convection equation in this case. However, if the velocity is
small (in our case, it is of order $\varepsilon^{2}$), then neglecting the
convection, in comparison with the diffusion term or iterating with respect to
the small term \footnote{The simplest coupled diffusion-convection equation
is
\[%
\begin{array}
[c]{l}%
-\mathrm{div}(\nu(C_{\varepsilon}){\mathcal{D}}u_{\varepsilon})+\nabla
p_{\varepsilon}=f(x)\\
\mathrm{div}u_{\varepsilon}=0,\\
-\Delta C_{\varepsilon}+ u_{\varepsilon}.\nabla C_{\varepsilon}=0
\end{array}
\]
where $C_{\varepsilon}$ is the concentration, $u_{\varepsilon}$ is the fluid
velocity and $p_{\varepsilon}$ is the pressure.\newline If $u_{\varepsilon}$
is small with respect to $C_{\varepsilon}$, we can write an expansion for
$C_{\varepsilon}$ with respect to the small parameter that is the ratio of
magnitudes of $|u_{\varepsilon}|$ and $|C_{\varepsilon}|$. Then for the terms
of this expansion, the third equation may be solved before the fluid motion
equations. Another possible approach is the successive approximations (fixed
point iterations) :
\[%
\begin{array}
[c]{l}%
-\mathrm{div}(\nu(C^{(n)}_{\varepsilon}){\mathcal{D}}u^{(n)}_{\varepsilon
})+\nabla p^{(n)}_{\varepsilon}=f(x)\\
\mathrm{div}u^{(n)}_{\varepsilon}=0,\\
-\Delta C^{(n)}_{\varepsilon}+ u^{(n-1)}_{\varepsilon}.\nabla C^{(n)}%
_{\varepsilon}=0
\end{array}
\]
where $n$ is the number of the iteration.\newline In both approaches, we get a
problem with the variable viscosity.}}, we get the steady state diffusion
equation; in absence of the source term in the right hand side, it has a
piecewise-linear asymptotic solution on the graph for the concentration. So,
in this simplified situation, the diffusion equation can be solved before the
fluid motion equation, and we obtain for the flow, the Stokes or Navier-Stokes
equation with a variable viscosity depending (via concentration) on the
longitudinal variable. These arguments provide the motivation for the Stokes
equation with variable viscosity; to our knowledge, it has not been studied
earlier from the asymptotic point of view.There are, of course, many other
practical problems involving fluids with variable viscosity. For example, the
presence of bacteria in suspension (see \cite{HABK}) may change locally the
viscosity.

In the first part, we consider the case of a flow in one rectangular channel
with the periodicity condition at the end of the channel. An asymptotic
expansion of solution is constructed and justified (a regular ansatz).

In the second part, we consider the case of a flow in one rectangular channel
with the inflow/outflow boundary conditions at the ends. In this case as well
we construct an asymptotic expansion of solution, that contains a regular
ansatz and the boundary layer correctors.

In the third part, we construct an asymptotic expansion for solution of the
Stokes equation set in a tube structure. Here as well we associate to each
rectangle a regular ansatz and then we glue all the regular ansatzes with help
of the boundary layer correctors exponential decaying from the junctions of
rectangles. This procedure is similar to the procedure described in
\cite{pana2005}.

In the fourth part we consider some numerical experiments comparing the
numerical and asymptotic solutions.

All constructed expansions are justified by calculation of residual terms and
application of \textit{a priori } estimates.

\section{Flow in one channel\label{sect2}}

Consider a small parameter $\varepsilon$, $\varepsilon=\frac{1}{q}$,
$q\in{\mathbb{N}}$, and define then a thin domain
\[
D_{\varepsilon}=\left\{  (x_{1},x_{2})\in{\mathbb{R}}^{2}:0<x_{1}%
<1,\ -\frac{\varepsilon}{2}<x_{2}<\frac{\varepsilon}{2}\right\}  .
\]
Assume that incompressible, viscous fluid fills the domain $D_{\varepsilon}$.
Let $f$ be the exterior force applied to the fluid.%

\begin{figure}
[h]
\begin{center}
\includegraphics[
height=1.0101in,
width=3.3512in
]%
{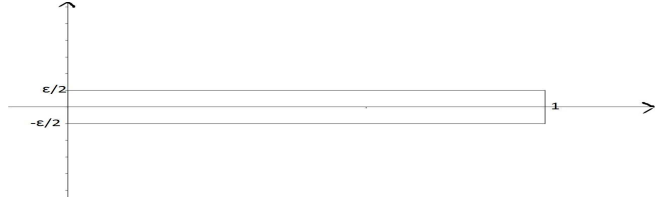}%
\caption{Thin domain}%
\label{f1}%
\end{center}
\end{figure}
%
Consider the following steady state Stokes problem:%
\begin{equation}
\left\{
\begin{array}
[c]{ll}%
-\mathrm{div}(\nu(x_{1}){\mathcal{D}}u_{\varepsilon})+\nabla p_{\varepsilon
}=f(x) & \text{in }D_{\varepsilon},\\
\mathrm{div}u_{\varepsilon}=0 & \text{in }D_{\varepsilon},\\
u_{\varepsilon}(x_{1},\frac{\varepsilon}{2})=0 & \text{for }x_{1}\in(0,1),\\
u_{\varepsilon}(x_{1},-\frac{\varepsilon}{2})=0 & \text{for }x_{1}\in(0,1),\\
u_{\varepsilon}(0,x_{2})=\varepsilon^{2}\varphi_{0}(\frac{x_{2}}{\varepsilon
}) & \text{for }x_{2}\in(-\frac{\varepsilon}{2},\frac{\varepsilon}{2}),\\
u_{\varepsilon}(1,x_{2})=\varepsilon^{2}\varphi_{1}(\frac{x_{2}}{\varepsilon
}) & \text{for }x_{2}\in(-\frac{\varepsilon}{2},\frac{\varepsilon}{2}),
\end{array}
\right.  \label{1.1}%
\end{equation}
The unknowns of this system are the velocity $u_{\varepsilon}$ and the
pressure $p_{\varepsilon}$ of the fluid.

The non homogeneous boundary conditions for the velocity are given by the
functions $\varphi_{0},\varphi_{1}$ whose second component is equal to zero,
and the first satisfies
\begin{equation}
\int_{-\frac{1}{2}}^{\frac{1}{2}}\varphi_{01}(\xi_{2})\mathrm{d}\xi_{2}%
=\int_{-\frac{1}{2}}^{\frac{1}{2}}\varphi_{11}(\xi_{2})\mathrm{d}\xi_{2},
\label{1.10}%
\end{equation}
where $\varphi_{i1}\in C_{0}^{2}([-\frac{1}{2},\frac{1}{2}])$. Here
$({\mathcal{D}}u_{\varepsilon})_{ij}=\frac{1}{2}\left(  \frac{\partial
u_{\varepsilon i}}{\partial x_{j}}+\frac{\partial u_{\varepsilon j}}{\partial
x_{i}}\right)  $, $\nu$ satisfies the following conditions
\begin{equation}
\nu(x_{1})=\nu_{0}+\nu_{1}(x_{1}) \label{1.11}%
\end{equation}
where $\nu_{1}\in C_{0}^{2}([0,1])$. Moreover there exist $\rho>0,$ $\kappa>0$
such that $\nu(x_{1})\geq\kappa$ for all $x_{1}\in(0,1)$ and $\nu_{1}%
(x_{1})=0$ for all $x_{1}\in(0,\rho)\cup(1-\rho,1)$ (and so, $\nu_{0}>0$).

\subsection{Variational formulation of the problem}

In order to obtain the variational formulation of problem \eqref{1.1}, let us
introduce the following space%
\begin{equation}
H(D_{\varepsilon})=\left\{  u\in(H_{0}^{1}(D_{\varepsilon}))^{2}%
:\mathrm{div}(u)=0\right\}  \label{1.12}%
\end{equation}
and assume that $f\in(L^{2}(D_{\varepsilon}))^{2}$.

Applying the extension theorem (see \cite{Girault-Raviart}), one can find a
function $\varphi^{\varepsilon}\in(H^{1}(D_{\varepsilon}))^{2}$ such that
$\mathrm{div}\varphi^{\varepsilon}=0$, $\varphi_{|_{x_{1}=0}}^{\varepsilon
}=\varphi_{0}$ and $\varphi_{|_{x_{1}=1}}^{\varepsilon}=\varphi_{1}$ and
$\varphi_{|_{x_{2}=\pm\frac{\varepsilon}{2}}}^{\varepsilon}=0$ . Changing the
unknown function $u_{\varepsilon}$ by $v_{\varepsilon}=u_{\varepsilon
}-\varepsilon^{2}\varphi^{\varepsilon},$ we give the variational formulation
of problem \eqref{1.1}:%
\begin{equation}
\int_{D_{\varepsilon}}\nu(x_{1}){\mathcal{D}}v_{\varepsilon}:{\mathcal{D}%
}{\psi}=\int_{D_{\varepsilon}}f\cdot\psi-\varepsilon^{2}\int_{D_{\varepsilon}%
}\nu(x_{1}){\mathcal{D}}\varphi_{\varepsilon}:{\mathcal{D}}{\psi}%
,\ \ \forall\psi\in H(D_{\varepsilon}). \label{2.1}%
\end{equation}

\begin{definition}
We say that $u_{\varepsilon}$ is a weak solution of problem \eqref{1.1} if
$v_{\varepsilon}=u_{\varepsilon}-\varepsilon^{2}\varphi^{\varepsilon}\in
H(D_{\varepsilon})$ and satisfies \eqref{2.1}.
\end{definition}

\begin{prop}
If $u_{\varepsilon}$ is a weak solution for problem \eqref{1.1}, then there
exists a distribution $p_{\varepsilon}\in{{\mathcal{D}}}^{^{\prime}%
}(D_{\varepsilon})$ such that $(u_{\varepsilon},p_{\varepsilon})$ satisfies
\eqref{1.1}$_{1}$ this problem in sense of distributions.

\end{prop}

\textit{Proof.} If we take $v_{\varepsilon}= u_{\varepsilon}- \varepsilon^{2}
\varphi^{\varepsilon}$, then from $(5)$ we have that
\[
\langle\mathrm{div}(\nu(x_{1}){\mathcal{D}} v_{\varepsilon}+ \varepsilon^{2}
\nu(x_{1}){\mathcal{D}} \varphi_{\varepsilon})+f , \psi\rangle= 0, \forall
\psi\in H(D_{\varepsilon})
\]
From De Rham lemma, it follows that there exists a distribution
$p_{\varepsilon}$, unique up to an additive constant, such that
\[
-\mathrm{div}(\nu(x_{1}){\mathcal{D}} v_{\varepsilon}) + \varepsilon
^{2}\mathrm{div}( \nu(x_{1}){\mathcal{D}} \varphi_{\varepsilon}) -f = - \nabla
p.
\]

\begin{theorem}
\label{th2.1}The variational problem \eqref{2.1} admits a unique solution
$v_{\varepsilon}\in H(D_{\varepsilon})$.
\end{theorem}

\textit{Proof.} The Riesz's theorem gives the existence and the uniqueness of
solution because the norms $\| v_{\varepsilon}\|_{I} = \sqrt
{\displaystyle{\int_{D_{\varepsilon}} \nu(x_{1}){\mathcal{D}} v_{\varepsilon
}:{\mathcal{D}} v_{\varepsilon}}}$ and $\| v_{\varepsilon}\|_{(H^{1}%
_{0}(D_{\varepsilon}))^{2}}$ are equivalent.

As a consequence, we have:

\begin{prop}
Let $v_{\varepsilon}\in H(D_{\varepsilon})$ be the solution of the variational
problem \eqref{2.1}. Then the following inequality holds%
\[
\Vert v_{\varepsilon}\Vert_{(H(D_{\varepsilon}))^{2}}\leq{\mathcal{C}}%
(\kappa,C^{\varepsilon}_{PF})\left(  \Vert f\Vert_{(L^{2}(D_{\varepsilon
}))^{2}}+\Vert\varphi^{\varepsilon}\Vert_{\left(  H^{1}(D_{\varepsilon
})\right)  ^{2}}\right)  ,
\]
{where $C^{\varepsilon}_{PF}$ stands for the Poincar\'{e}-Friedrichs
inequality constant and $\kappa$ is the lower bound of the viscosity
\eqref{1.11}}.
\end{prop}

Mention that $C^{\varepsilon}_{PF}$ can be estimated by $\varepsilon
\tilde{{\mathcal{C}}}$, where $\tilde{{\mathcal{C}}}$ is independent of
$\varepsilon$ (see \cite{pana2005}).

Consider the case when $\varphi^{\varepsilon}=0.$ Then we have:

\begin{prop}
The following inequality holds%
\[
\Vert\nabla p_{\varepsilon}\Vert_{H^{-1}(D_{\varepsilon})}\leq C\Vert
f\Vert_{\left(  L^{2}(D_{\varepsilon})\right)  ^{2}},
\]
where $C$ is a constant independent of $\varepsilon.$
\end{prop}

\textit{Proof.} From De Rham lemma, it follows that there exists a
distribution $p_{\varepsilon}$, unique up to an additive constant, such that
\[
\nabla p_{\varepsilon}=\mathrm{div}(\nu(x_{1}){\mathcal{D}} u_{\varepsilon}) +
f.
\]
It means that
\[
\int_{D_{\varepsilon}} p_{\varepsilon} \mathrm{div} \psi= \int_{D_{\varepsilon
}} \nu(x_{1}){\mathcal{D}} u_{\varepsilon}:{\mathcal{D}} {\psi} -
\int_{D_{\varepsilon}}f\cdot\psi, \quad\forall\psi\in H(D_{\varepsilon}))^{2},
\]
and so,
\[
\sup_{\psi\in(H^{1}_{0}(D_{\varepsilon}))^{2}} \frac{\vert\int_{D_{\varepsilon
}} p_{\varepsilon} div \psi\vert} { \Vert\nabla\psi\Vert_{(L^{2}%
(D_{\varepsilon}))^{4}}} = \sup_{\psi\in(H^{1}_{0}(D_{\varepsilon}))^{2}}
\frac{\vert\int_{D_{\varepsilon}} \nu(x_{1}){\mathcal{D}} u_{\varepsilon
}:{\mathcal{D}} {\psi} - \int_{D_{\varepsilon}}f\cdot\psi\vert}{ \Vert
\nabla\psi\Vert_{(L^{2}(D_{\varepsilon}))^{4}}} \leq
\]
\[
\leq({\mathcal{C}}(\kappa,C^{\varepsilon}_{PF})+C^{\varepsilon}_{PF}) \| f
\|_{(L^{2}(D_{\varepsilon}))^{2}},
\]
where ${\mathcal{C}}(\kappa,C^{\varepsilon}_{PF}) $ is the constant of the a
priori estimate for $u_{\varepsilon} $ (see above) and $C^{\varepsilon}_{PF}$
is the Poincar\'e-Friedrichs constant (it is of order $\varepsilon$). This
inequality gives the estimate of $\nabla p_{\varepsilon} $ in the $H^{-1}%
-$norm.

An asymptotic analysis of problem \eqref{1.1} shows that an asymptotic
solution is given by a Poiseuille type flow, with two boundary layer
correctors localized in some neighborhoods of the ends of the channel.
Didactically, it would be better to separate the construction of the
Poiseuille type flow for varying viscosity and the construction of the
boundary layer correctors. That's why we simplify problem $\eqref{1.1}$, and
replace $\eqref{1.1}_{5}$, $\eqref{1.1}_{6}$ by the periodicity condition with
respect to $x_{1}$.

Introduce the Sobolev space
\[
H_{per}(D_{\varepsilon})=\left\{  u\in{\left(  H_{per,1,0}^{1}(D_{\varepsilon
})\right)  }^{2}:\operatorname{div}u=0\right\}  .
\]
Here $H_{per,1,0}^{1}(D_{\varepsilon})$ is the completion (with respect to the
$H^{1}(D_{\varepsilon})$-norm) of the space of the $C^{\infty}({\mathbb{R}%
}\times\lbrack-\frac{\varepsilon}{2},\frac{\varepsilon}{2}])$-functions
vanishing at the boundary $x_{2}=\pm\frac{\varepsilon}{2}$ and 1-periodic in
$x_{1}$. As in the beginning of the section, $f\in(L^{2}(D_{\varepsilon}%
))^{2}$.

\begin{definition}
We say that $u_{\varepsilon}\in H_{per}(D_{\varepsilon})$ is a weak solution
of the periodic problem%

\begin{equation}
\left\{
\begin{array}
[c]{ll}%
-\mathrm{\operatorname{div}}(\nu(x_{1}){\mathcal{D}}u_{\varepsilon})+\nabla
p_{\varepsilon}=f(x), & \text{in }D_{\varepsilon},\\
\mathrm{\operatorname{div}}u_{\varepsilon}=0, & \text{in }D_{\varepsilon},\\
u_{\varepsilon}(x_{1},\frac{\varepsilon}{2})=0, & \text{for }x_{1}\in(0,1),\\
u_{\varepsilon}(x_{1},-\frac{\varepsilon}{2})=0, & \text{for }x_{1}\in(0,1),\\
u_{\varepsilon}\text{ is }1-\text{periodic in }x_{1}, &
\end{array}
\right.  \label{A.1}%
\end{equation}
if and only if it satisfies the integral identity
\begin{equation}
\int_{D_{\varepsilon}}\nu(x_{1}){\mathcal{D}}u_{\varepsilon}:{\mathcal{D}%
}{\psi}=\int_{D_{\varepsilon}}f\cdot\psi,\text{ \ }\forall\psi\in
H_{per}(D_{\varepsilon}). \label{A.11}%
\end{equation}

\end{definition}

As in theorem \ref{th2.1}, we apply the Riesz theorem and prove the existence
and the uniqueness of $u_{\varepsilon}\in H_{per}(D_{\varepsilon})$, solution
of problem \eqref{A.1}. The \textit{a priori} estimate is given by
\[
\Vert u_{\varepsilon}\Vert_{(H^{1}(D_{\varepsilon}))^{2}}\leq{\mathcal{C}%
}(\kappa,C^{\varepsilon}_{PF})\Vert f\Vert_{(L^{2}(D_{\varepsilon}))^{2}}%
\]
{where ${\mathcal{C}}(\kappa,C^{\varepsilon}_{PF})=O\left(  \varepsilon\right)
$ and $\kappa$ is the lower bound of the viscosity \eqref{1.11}}.

\begin{prop}
If $u_{\varepsilon}$ is a weak solution for problem \eqref{A.1} then there
exists a distribution $p_{\varepsilon}\in{{\mathcal{D}}}^{^{\prime}%
}(D_{\varepsilon})$ such that $(u_{\varepsilon},p_{\varepsilon})$ satisfies
\eqref{A.1} $_{1}$  in sense of distributions, and
\[
\Vert\nabla p_{\varepsilon}\Vert_{H^{-1}(D_{\varepsilon})}\leq C\Vert
f\Vert_{\left(  L^{2}(D_{\varepsilon})\right)  ^{2}},
\]
where $C$ is a constant independent of $\varepsilon.$
\end{prop}

\textit{Proof.} By $\eqref{A.11}$
\[
\int_{D_{\varepsilon}} \nu(x_{1}){\mathcal{D}} u_{\varepsilon}:{\mathcal{D}}
{\psi} = \int_{D_{\varepsilon}}f\cdot\psi, \qquad\forall\psi\in H_{per}%
(D_{\varepsilon}).
\]
From De Rham lemma, it follows that there exists a distribution
$p_{\varepsilon}$, unique up to an additive constant, such that
\[
\nabla p_{\varepsilon}=\mathrm{div}(\nu(x_{1}){\mathcal{D}} u_{\varepsilon}) +
f.
\]
It means that
\[
\int_{D_{\varepsilon}} p_{\varepsilon} \mathrm{div} \psi= \int_{D_{\varepsilon
}} \nu(x_{1}){\mathcal{D}} u_{\varepsilon}:{\mathcal{D}} {\psi} -
\int_{D_{\varepsilon}}f\cdot\psi\qquad\forall\psi\in(H^{1}_{per,1,0}%
(D_{\varepsilon}))^{2},
\]
and so,
\[
\sup_{\psi\in(H^{1}_{per,1,0}(D_{\varepsilon}))^{2}} \frac{\vert
\int_{D_{\varepsilon}} p_{\varepsilon} \mathrm{div} \psi\vert}{ \Vert
\nabla\psi\Vert_{(L^{2}(D_{\varepsilon}))^{4}}} = \sup_{\psi\in(H^{1}%
_{per,1,0}(D_{\varepsilon}))^{2}} \frac{\vert\int_{D_{\varepsilon}} \nu
(x_{1}){\mathcal{D}} u_{\varepsilon}:{\mathcal{D}} {\psi} - \int
_{D_{\varepsilon}}f.\psi\vert}{ \Vert\nabla\psi\Vert_{(L^{2}(D_{\varepsilon
}))^{4}}} \leq
\]
\[
\leq({\mathcal{C}}(\kappa,C^{\varepsilon}_{PF})+C^{\varepsilon}_{PF}) \| f
\|_{(L^{2}(D_{\varepsilon}))^{2}},
\]
where ${\mathcal{C}}(\kappa,C^{\varepsilon}_{PF}) $ is the constant of the a
priori estimate for $u_{\varepsilon} $ (see above) and $C^{\varepsilon}_{PF}$
is the Poincar\'e-Friedrichs constant (it is of order $\varepsilon$). This
inequality gives the estimate of $\nabla p_{\varepsilon} $ in the $H^{-1}%
-$norm.

\subsection{Asymptotic analysis of the problem\label{sect2.1}}

Let us first construct the asymptotic expansion for the solution of the
periodic problem \eqref{A.1}; then we will study the non periodic one.

Define the infinite layer
\[
\Omega_{\varepsilon}=\left\{  (x_{1},x_{2})\in{\mathbb{R}}^{2}:x_{1}%
\in{\mathbb{R}},-\frac{\varepsilon}{2}<x_{2}<\frac{\varepsilon}{2}\right\}  .
\]
Assume that $f=f_{1}e_{1}$, $f_{1}\in C^{\infty}({\mathbb{R}})$, $f_{1}$ is
1-periodic in $x_{1}$.%
\begin{figure}
[ptb]
\begin{center}
\includegraphics[
height=1.1243in,
width=3.3537in
]%
{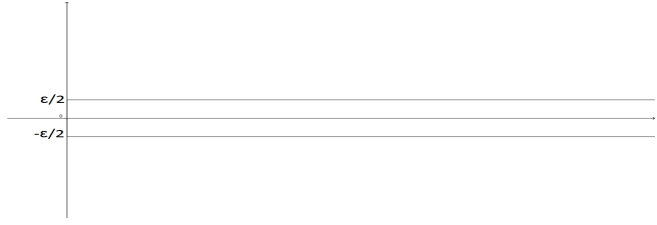}%
\caption{Infinite layer}%
\label{f2}%
\end{center}
\end{figure}

Denote ${\langle\psi\rangle}_{1}=\int_{0}^{1}\psi(x_{1},x_{2})\mathrm{d}x_{1}%
$, $\langle\psi{\rangle}_{2}=\int_{-\frac{1}{2}}^{\frac{1}{2}}\psi(x_{1}%
,x_{2})\mathrm{d}x_{2}$. An asymptotic solution is written as:%
\begin{equation}
\left\{
\begin{array}
[c]{lcl}%
u_{1}^{k}(x_{1},x_{2}) & = & \displaystyle{\sum_{j=0}^{k}}\varepsilon
^{j+2}u_{1,j}\left(  x_{1},\frac{x_{2}}{\varepsilon}\right)  ,\\
u_{2}^{k}(x_{1},x_{2}) & = & \displaystyle{\sum_{j=0}^{k}}\varepsilon
^{j+3}u_{2,j}\left(  x_{1},\frac{x_{2}}{\varepsilon}\right)  ,\\
p^{k}(x_{1},x_{2}) & = & \displaystyle{\sum_{j=0}^{k}}\varepsilon^{j+1}%
p_{j}\left(  x_{1},\frac{x_{2}}{\varepsilon}\right)  +\displaystyle{\sum
_{j=0}^{k}}\varepsilon^{j}q_{j}(x_{1}),
\end{array}
\right.  \label{A.2}%
\end{equation}
with $u_{1,j},u_{2,j},p_{j}$ and $q_{j}$ 1-periodic functions in $x_{1}$, such
that $\langle p_{j}\rangle_{2}=0$.

Substituting \eqref{A.2}, in \eqref{A.1}, equating the coefficients of the
same powers of $\varepsilon$ and denoting $\xi_{2}=\frac{x_{2}}{\varepsilon}$
we obtain:%
\begin{equation}
\left\{
\begin{array}
[c]{l}%
-\frac{\partial}{\partial x_{1}}\left(  \nu(x_{1})\frac{\partial u_{1,j-2}%
}{\partial x_{1}}\right)  -\frac{\nu(x_{1})}{2}\frac{\partial^{2}u_{1,j}%
}{\partial\xi_{2}^{2}}-\frac{\nu(x_{1})}{2}\frac{\partial^{2}u_{2,j-2}%
}{\partial\xi_{2}\partial x_{1}}+\frac{\partial p_{j-1}}{\partial x_{1}}%
+\frac{\partial q_{j}}{\partial x_{1}}=f_{1}\delta_{j0},\\
\\
\frac{1}{2}\frac{\partial}{\partial x_{1}}\left[  \nu(x_{1})\left(
\frac{\partial u_{1,j-1}}{\partial\xi_{2}}+\frac{\partial u_{2,j-3}}{\partial
x_{1}}\right)  \right]  +\nu(x_{1})\frac{\partial^{2}u_{2,j-1}}{\partial
\xi_{2}^{2}}-\frac{\partial p_{j}}{\partial\xi_{2}}=0,\\
\\
\frac{\partial u_{1,j}}{\partial x_{1}}+\frac{\partial u_{2,j}}{\partial
\xi_{2}}=0,\\
\\
u_{1,j}\left(  x_{1},\pm\frac{1}{2}\right)  =0,\\
\\
u_{2,j}\left(  x_{1},\pm\frac{1}{2}\right)  =0.
\end{array}
\right.  \label{A.3}%
\end{equation}
Denote:
\[
N_{1}(\xi_{2})=\frac{1}{2}\left(  \xi_{2}^{2}-\frac{1}{4}\right)  ,
\]
(it satisfies $N_{1}^{^{\prime\prime}}=1$, $N_{1}(\pm\frac{1}{2})=0$) and
\[
N_{2}(\xi_{2})=\int_{-\frac{1}{2}}^{\xi_{2}}N_{1}(\tau)\mathrm{d}\tau,
\]
(here $N_{2}(\frac{1}{2})={%
{\displaystyle\int_{-\frac{1}{2}}^{\frac{1}{2}}}
N_{1}(\tau)\mathrm{d}\tau}=-\frac{1}{12}$). Denote:
\[
\left\{
\begin{array}
[c]{l}%
D^{-1}:F\longrightarrow\displaystyle{\int_{-\frac{1}{2}}^{\xi_{2}}F(x_{1}%
,\tau)\mathrm{d}\tau,}\\
\tilde{D}^{-1}:F\longrightarrow\displaystyle{\int_{-\frac{1}{2}}^{\xi_{2}%
}F(x_{1},\tau)\mathrm{d}\tau-\int_{-\frac{1}{2}}^{\frac{1}{2}}\int_{-\frac
{1}{2}}^{\theta}F(x_{1},\tau)\mathrm{d}\tau\mathrm{d}\theta,}\\
D^{-2}:F\longrightarrow\displaystyle{\int_{-\frac{1}{2}}^{\xi_{2}}\int
_{-\frac{1}{2}}^{\theta}F(x_{1},\tau)\mathrm{d}\tau\mathrm{d}\theta-\left(
\xi_{2}+\frac{1}{2}\right)  \int_{-\frac{1}{2}}^{\frac{1}{2}}\int_{-\frac
{1}{2}}^{\theta}F(x_{1},\tau)\mathrm{d}\tau\mathrm{d}\theta.}%
\end{array}
\right.
\]

\begin{theorem}
The unknowns of \eqref{A.3} $u_{1,j},u_{2,j},p_{j},q_{j}$ are given by the
following relations:
\begin{equation}
\left\{
\begin{array}
[c]{lcl}%
u_{1,j} & = & -D^{-2}\left\{  \frac{\partial^{2}u_{2,j-2}}{\partial\xi
_{2}\partial x_{1}}+\frac{2}{\nu(x_{1})}\left(  \frac{\partial}{\partial
x_{1}}\left(  \nu(x_{1})\frac{\partial u_{1,j-2}}{\partial x_{1}}\right)
-\frac{\partial p_{j-1}}{\partial x_{1}}\right)  \right\} \\
&  & +\frac{2}{\nu(x_{1})}N_{1}(\xi_{2})\left(  \frac{\partial q_{j}}{\partial
x_{1}}-f_{1}\delta_{j0}\right)  ,\\
&  & \\
u_{2,j} & = & D^{-1}D^{-2}\left\{  \frac{\partial^{3}u_{2,j-2}}{\partial
\xi_{2}\partial x_{1}^{2}}+\frac{\partial}{\partial x_{1}}\left(  \frac{2}%
{\nu(x_{1})}\left(  \frac{\partial}{\partial x_{1}}\left(  \nu(x_{1}%
)\frac{\partial u_{1,j-2}}{\partial x_{1}}\right)  -\frac{\partial p_{j-1}%
}{\partial x_{1}}\right)  \right)  \right\} \\
&  & -\frac{\partial}{\partial x_{1}}\left(  \frac{2}{\nu(x_{1})}\left(
\frac{\partial q_{j}}{\partial x_{1}}-f_{1}\delta_{j0}\right)  \right)
N_{2}(\xi_{2}),\\
&  & \\
p_{j} & = & \tilde{D}^{-1}\left\{  \frac{1}{2}\frac{\partial}{\partial x_{1}%
}\left[  \nu(x_{1})\left(  \frac{\partial u_{1,j-1}}{\partial\xi_{2}}%
+\frac{\partial u_{2,j-3}}{\partial x_{1}}\right)  \right]  +\nu(x_{1}%
)\frac{\partial^{2}u_{2,j-1}}{\partial\xi_{2}^{2}}\right\}  .
\end{array}
\right.  \label{A.4}%
\end{equation}

\end{theorem}

\textit{Proof.} Integrating twice $\eqref{A.3}_{1}$ and using boundary
conditions $\eqref{A.3}_{4}$, we get $\eqref{A.4}_{1}$. This relation gives an
expression of the unknown $u_{1,j}$ via $q_{j}$. All other functions contained
by this relation are either known from previous computations or equal to zero.
We integrate next the incompressibility condition $\eqref{A.3}_{3}$ with
respect to $\xi_{2}$ with the boundary condition$\eqref{A.3}_{5}$ and get
$u_{2,j}$ . Finally, integrating $\eqref{A.3}_{2}$ we get $p_{j}$. The unknown
function $q_{j}$ is determined from the boundary condition $u_{2,j}\left(
x_{1},\frac{1}{2}\right)  =0$. Actually, the boundary conditions
$\eqref{A.3}_{4}$ and $u_{2,j}\left(  x_{1},-\frac{1}{2}\right)  =0$ are
satisfied by the definition of $D^{-2}$ and $D^{-1}$, while $u_{2,j}\left(
x_{1},\frac{1}{2}\right)  =0$ gives:
\begin{align}
-\frac{\partial}{\partial x_{1}}  &  \left(  \frac{1}{6\nu(x_{1})}\left(
\frac{\partial q_{j}}{\partial x_{1}}-f_{1}\delta_{j0}\right)  \right)
=\label{A.5}\\
&  ={D^{-1}D^{-2}\left\{  \frac{\partial^{3}u_{2,j-2}}{\partial\xi_{2}\partial
x_{1}^{2}}+\frac{\partial}{\partial x_{1}}\left\{  \frac{2}{\nu(x_{1})}\left(
\frac{\partial}{\partial x_{1}}\left(  \nu(x_{1})\frac{\partial u_{1,j-2}%
}{\partial x_{1}}\right)  -\frac{\partial p_{j-1}}{\partial x_{1}}\right)
\right\}  \right\}  }_{\left\vert _{\xi_{2}=\frac{1}{2}}\right.  . }\nonumber
\end{align}

In particular, for $j=0$, we have the Darcy equation for the leading term of
the pressure $q_{0}$:%
\begin{equation}
\frac{\partial}{\partial x_{1}}\left(  \frac{1}{6\nu(x_{1})}\left(
\frac{\partial q_{0}}{\partial x_{1}}-f_{1}\right)  \right)  =0, \label{A.6}%
\end{equation}
and so $q_{0}$ is the 1-periodic function given by the formula
\[
q_{0}(x_{1}) = \displaystyle{\int_{0}^{x_{1}}\left( f_{1}(\tau)-\frac{\langle
f_{1}\rangle_{1}}{\langle\nu\rangle_{1}} \nu(\tau)\right) \mathrm{d}\tau-
\left\langle \int_{0}^{x_{1}}\left( f_{1}(\tau)-\frac{\langle f_{1}
\rangle_{1}}{\langle\nu\rangle_{1}}\nu(\tau)\right) \mathrm{d}\tau
\right\rangle _{1}}.
\]
Then $u_{1,0}(x_{1},\xi_{2}) = \frac{2}{\nu(x_{1})}\left( \frac{\partial
q_{0}}{\partial x_{1}} - f_{1}\right) N_{1}(\xi_{2})$, $u_{2,0}(x_{1},\xi_{2})
= 0$ and $p_{0}(x_{1},\xi_{2}) = 0$.\newline For $j=1$ $q_{1}$ is the
1-periodic solution of equation
\[
\frac{\partial}{\partial x_{1}}\left(  \frac{1}{6\nu(x_{1})} \frac{\partial
q_{1}}{\partial x_{1}} \right)  = 0.
\]
It follows: $q_{1}=0$, $u_{1,1}(x_{1},\xi_{2})=0 \mbox{ and } u_{2,1}%
(x_{1},\xi_{2}) = 0$. For $j=1$, $\eqref{A.4}_{3}$ yields:
\[
p_{1}(x_{1},\xi_{2}) = \frac{\partial}{\partial x_{1}} \left(  \frac{\partial
q_{0}}{\partial x_{1}} - f_{1}\right) \left(  N_{1}(\xi_{2}) - \langle
N_{1}\rangle_{2}\right) .
\]
For $j=2$, \eqref{A.5} becomes
\begin{equation}
\frac{\partial}{\partial x_{1}}\left(  \frac{1}{6\nu(x_{1})} \frac{\partial
q_{2}}{\partial x_{1}}\right)  = \frac{\partial}{\partial x_{1}}\left(
\frac{2}{\nu(x_{1})} \frac{\partial}{\partial x_{1}}\left(  \frac{\partial
q_{0}} {\partial x_{1}} - f_{1} \right) \right) D^{-1}D^{-2}\left(  N_{1}%
(\xi_{2}) -\langle N_{1}\rangle_{2}\right)  _{\left| _{\xi_{2} = \frac{1}{2}}
\right. }, \label{A.70}%
\end{equation}
where $D^{-1}D^{-2}\left(  N_{1}(\xi_{2}) -\langle N_{1}\rangle_{2}\right)
_{\left| _{\xi_{2} = \frac{1}{2}} \right. } = {\mathcal{R}} \neq0$, $q_{2}$ is
1-periodic solution of given by:
\[
q_{2}(x_{1})=-12{\mathcal{R}}{\langle f_{1} \rangle_{1}} \left( \frac
{\nu(x_{1})}{\langle\nu\rangle_{1}} -1\right) .
\]
It follows:
\[
\left\{
\begin{array}
[c]{l}%
u_{1,2}(x_{1},\xi_{2}) = \frac{2}{\nu(x_{1})} \frac{\partial^{2}}{\partial
x_{1}^{2}} (\frac{\partial q_{0}}{\partial x_{1}} - f_{1}) D^{-2}\left(
N_{1}(\frac{x_{2}}{\varepsilon}) -\langle N_{1}\rangle_{2} \right)  +\frac
{2}{\nu(x_{1})} N_{1}(\frac{x_{2}} {\varepsilon}) \frac{\partial q_{2}%
}{\partial x_{1}},\\
\\
u_{2,2}(x_{1},\xi_{2})=- \frac{\partial}{\partial x_{1}}\left[  \frac{2}%
{\nu(x_{1})}\frac{\partial^{2}}{\partial x_{1}^{2}} \left(  \frac{\partial
q_{0}}{\partial x_{1}} - f_{1}\right) \right] D^{-1}D^{-2}\left(  N_{1}%
(\frac{x_{2}}{\varepsilon})-\langle N_{1}\rangle_{2}\right)  - \frac{\partial
}{\partial x_{1}}\left( \frac{2}{\nu(x_{1})} N_{2}(\frac{x_{2}}{\varepsilon})
\frac{\partial q_{2}}{\partial x_{1}}\right) ,
\end{array}
\right.
\]
and we get the following theorem.

\begin{theorem}
The asymptotic solution of problem \eqref{A.1} is given by:%
\begin{equation}
\left\{
\begin{array}
[c]{lcl}%
u_{1}^{k}(x_{1},x_{2}) & = & \varepsilon^{2}\frac{2}{\nu(x_{1})}\left(
\frac{\partial q_{0}}{\partial x_{1}}-f_{1}\right)  N_{1}(\frac{x_{2}%
}{\varepsilon})\\
&  & +\varepsilon^{4}\left(  \frac{2}{\nu(x_{1})}\frac{\partial^{2}}{\partial
x_{1}^{2}}(\frac{\partial q_{0}}{\partial x_{1}}-f_{1})D^{-2}\left(
N_{1}(\frac{x_{2}}{\varepsilon})-\langle N_{1}\rangle_{2}\right)  +\frac
{2}{\nu(x_{1})}N_{1}(\frac{x_{2}}{\varepsilon})\frac{\partial q_{2}}{\partial
x_{1}}\right)  +O(\varepsilon^{5}),\\
&  & \\
u_{2}^{k}(x_{1},x_{2}) & = & -\varepsilon^{5}\left(  \frac{\partial}{\partial
x_{1}}\left[  \frac{2}{\nu(x_{1})}\frac{\partial^{2}}{\partial x_{1}^{2}%
}\left(  \frac{\partial q_{0}}{\partial x_{1}}-f_{1}\right)  \right]
D^{-1}D^{-2}\left(  N_{1}(\frac{x_{2}}{\varepsilon})-\langle N_{1}\rangle
_{2}\right)  \right. \\
&  & \left.  +\frac{\partial}{\partial x_{1}}\left(  \frac{2}{\nu(x_{1})}%
N_{2}(\frac{x_{2}}{\varepsilon})\frac{\partial q_{2}}{\partial x_{1}}\right)
\right)  +O(\varepsilon^{6}),\\
&  & \\
p^{k}(x_{1},x_{2}) & = & q_{0}+\varepsilon^{2}\left(  q_{2}(x_{1}%
)+\frac{\partial}{\partial x_{1}}\left(  \frac{\partial q_{0}}{\partial x_{1}%
}-f_{1}\right)  \left(  N_{1}(\frac{x_{2}}{\varepsilon})-\langle N_{1}%
\rangle_{2}\right)  \right)  +O(\varepsilon^{3}),
\end{array}
\right.  \label{A.7}%
\end{equation}
where $q_{0}$ and $q_{2}$ are periodic solutions of obtained as the unique
solution of \eqref{A.5} \eqref{A.70} respectively.
\end{theorem}

\begin{rque}
Solution given by \eqref{A.7} can be written as follows:
\[
\left\{
\begin{array}
[c]{lcl}%
u_{1}^{k}(x_{1},x_{2}) & = & -\varepsilon^{2}\frac{2\langle f_{1}\rangle_{1}%
}{\langle\nu(x_{1})\rangle_{1}}N_{1}(\xi_{2})-\varepsilon^{4}\left(  \frac
{2}{\nu(x_{1})}\frac{\langle f_{1}\rangle_{1}}{\langle\nu\rangle_{1}}%
\nu^{^{\prime\prime}}(x_{1})D^{-2}\left(  N_{1}(\frac{x_{2}}{\varepsilon
})-\langle N_{1}\rangle_{2}\right)  \right. \\
&  & \left.  +\frac{24{\mathcal{R}}\langle f_{1}\rangle_{1}}{\langle\nu
(x_{1})\rangle_{1}}\frac{\nu^{^{\prime}}(x_{1})}{\nu(x_{1})}N_{1}(\frac{x_{2}%
}{\varepsilon})\right)  +O(\varepsilon^{5}),\\
&  & \\
u_{2}^{k}(x_{1},x_{2}) & = & \varepsilon^{5}\left(  \frac{\partial}{\partial
x_{1}}\left[  \frac{2}{\nu(x_{1})}\frac{\langle f_{1}\rangle_{1}}{\langle
\nu\rangle_{1}}\nu^{^{\prime\prime}}(x_{1})\right]  D^{-1}D^{-2}\left(
N_{1}(\frac{x_{2}}{\varepsilon})-\langle N_{1}\rangle_{2}\right)  \right. \\
&  & \left.  +\frac{24{\mathcal{R}}\langle f_{1}\rangle_{1}}{\langle\nu
(x_{1})\rangle_{1}}(\frac{\nu^{^{\prime}}(x_{1})}{\nu(x_{1})})^{\prime}%
N_{2}(\frac{x_{2}}{\varepsilon})\right)  +O(\varepsilon^{6}),\\
&  & \\
p^{k}(x_{1},x_{2}) & = & \displaystyle{\int_{0}^{x_{1}}\left(  f_{1}%
(\tau)-\frac{\langle f_{1}\rangle_{1}}{\langle\nu\rangle_{1}}\nu(\tau)\right)
\mathrm{d}\tau-\left\langle \int_{0}^{x_{1}}\left(  f_{1}(\tau)-\frac{\langle
f_{1}\rangle_{1}}{\langle\nu\rangle_{1}}\nu(\tau)\right)  \mathrm{d}%
\tau\right\rangle _{1}}\\
&  & +\varepsilon^{2}\left(  -12{\mathcal{R}}{\langle f_{1}\rangle_{1}}\left(
\frac{\nu(x_{1})}{\langle\nu\rangle_{1}}-1\right)  -\frac{\langle f_{1}%
\rangle_{1}}{\langle\nu\rangle_{1}}\nu^{^{\prime}}(x_{1})\left(  N_{1}%
(\frac{x_{2}}{\varepsilon})-\langle N_{1}\rangle_{2}\right)  \right)
+O(\varepsilon^{3}).
\end{array}
\right.
\]

\end{rque}

Introduce the following function:%
\begin{equation}%
\begin{array}
[c]{lcl}%
F^{k}(x_{1},x_{2}) & = & \left(  \frac{\partial}{\partial x_{1}}\left(
\nu(x_{1})\frac{\partial u_{1,k-1}}{\partial x_{1}}\right)  +\frac{\nu(x_{1}%
)}{2}\frac{\partial^{2}u_{2,k-1}}{\partial x_{1}\partial\xi_{2}}%
-\frac{\partial p_{k}}{\partial x_{1}}\right)  e_{1}\\
&  & \\
& + & \left(  \frac{1}{2}\frac{\partial}{\partial x_{1}}\left[  \nu
(x_{1})\left(  \frac{\partial u_{1,k}}{\partial\xi_{2}}+\frac{\partial
u_{2,k-2}}{\partial x_{1}}\right)  \right]  +\nu(x_{1})\frac{\partial
^{2}u_{2,k}}{\partial\xi_{2}^{2}}\right)  e_{2}\\
&  & \\
& + & \varepsilon\left(  \frac{\partial}{\partial x_{1}}\left(  \nu
(x_{1})\frac{\partial u_{1,k}}{\partial x_{1}}\right)  +\frac{\nu(x_{1})}%
{2}\frac{\partial^{2}u_{2,k}}{\partial x_{1}\partial\xi_{2}}\right)
e_{1}+\left(  \frac{1}{2}\frac{\partial}{\partial x_{1}}\left(  \nu
(x_{1})\frac{\partial u_{2,k-1}}{\partial x_{1}}\right)  \right)  e_{2}\\
&  & \\
& + & \varepsilon^{2}\left(  \frac{1}{2}\frac{\partial}{\partial x_{1}}\left(
\nu(x_{1})\frac{\partial u_{2,k}}{\partial x_{1}}\right)  \right)  e_{2}.
\end{array}
\label{rol}%
\end{equation}

\begin{theorem}
Let $(u^{k},p^{k})$ be the asymptotic solution given by \eqref{A.2} and
$(u,p)$ the solution of \eqref{A.1}. Then the following estimate holds:
\[
\Vert u^{k}-u\Vert_{{H^{1}(D_{\varepsilon})}^{2}}=O\left(  \varepsilon
^{k+\frac{5}{2}}\right)  .
\]

\end{theorem}

\textit{Proof.} Denote $(U^{k},P^{k}) = (u^{k} - u, p^{k} -p)$. We obtain the
following problems for $(U^{k},P^{k})$:
\begin{equation}
\label{A.8}\left\{
\begin{array}
[c]{ll}%
-\mathrm{div}(\nu(x_{1}){\mathcal{D}} U^{k}) + \nabla P^{k}= -\varepsilon
^{k+1}F^{k} (x_{1},x_{2}) & \quad in \quad D_{\varepsilon},\\
\mathrm{div} U^{k}= 0 & \quad in \quad D_{\varepsilon},\\
U^{k}(x_{1},\frac{\varepsilon}{2})=0 & \quad for \quad x_{1}\in(0,1),\\
U^{k}(x_{1},-\frac{\varepsilon}{2})= 0 & \quad for \quad x_{1}\in(0,1),\\
U^{k} \quad is \quad1-periodic \quad in \quad x_{1}. &
\end{array}
\right.
\end{equation}
The a priori estimates give us:
\[
\| u^{k}-u\|_{{H^{1}(D_{\varepsilon})}^{2}} \leq{\mathcal{C}}(K,C_{PF}%
)\|\varepsilon^{k+1}F^{k}\|_{{L^{2}(D_{\varepsilon})}^{2}} = O\left(
\varepsilon^{k+\frac{5}{2}} \right)
\]
and
\[
\|\nabla p^{k} -\nabla p\|_{H^{-1}(D_{\varepsilon})}=O\left( \varepsilon
^{k+\frac{5}{2}} \right) .
\]

Let us consider now the non periodic case, i.e. let us construct an asymptotic
expansion of the solution for problem \eqref{1.1}. Assume
\[
f=f_{1}(x_{1})e_{1},\ \ f_{1}\in C^{\infty}([0,1]).
\]

Define an asymptotic solution by:%
\begin{equation}
\left\{
\begin{array}
[c]{lcl}%
\hat{u}^{k}(x_{1},x_{2}) & = & {u}^{k}(x_{1},x_{2})+{u}_{BL0}^{k}(\frac
{x}{\varepsilon})+{u}_{BL1}^{k}(\frac{x_{1}-1}{\varepsilon},\frac{x_{2}%
}{\varepsilon})\\
\hat{p}^{k}(x_{1},x_{2}) & = & {p}^{k}(x_{1},x_{2})+{p}_{BL0}^{k}(\frac
{x}{\varepsilon})+{p}_{BL1}^{k}(\frac{x_{1}-1}{\varepsilon},\frac{x_{2}%
}{\varepsilon})
\end{array}
\right.  \label{B.1}%
\end{equation}
The expressions of ${u}^{k}$, ${p}^{k}$ are the same as above: \eqref{A.2}
,\eqref{A.4} ,\eqref{A.5}

Since the functions given by \eqref{A.2} do not satisfy exactly the boundary
conditions $\eqref{1.1}_{5},\eqref{1.1}_{6}$ we add the boundary layer
correctors. They correspond to the left end for $i=0$ and to the right end for
$i=1$ and their expressions are given by:%
\begin{equation}
\left\{
\begin{array}
[c]{lcl}%
{u}_{BLi}^{k}(\frac{x_{1}-i}{\varepsilon},\frac{x_{2}}{\varepsilon}) & = &
\displaystyle{\sum_{j=0}^{k}}\varepsilon^{j+2}u_{j}^{i}\left(  \frac{x_{1}%
-i}{\varepsilon},\frac{x_{2}}{\varepsilon}\right) \\
{p}_{BLi}^{k}(\frac{x_{1}-i}{\varepsilon},\frac{x_{2}}{\varepsilon}) & = &
\displaystyle{\sum_{j=0}^{k}}\varepsilon^{j+1}p_{j}^{i}\left(  \frac{x_{1}%
-i}{\varepsilon},\frac{x_{2}}{\varepsilon}\right)
\end{array}
\right.  \label{B.3}%
\end{equation}
Here $i=0,1.$ To obtain problems for the boundary layers corresponding to the
left side, we introduce the domain $\Pi^{+}=(0,\infty)\times(-\frac{1}%
{2},\frac{1}{2})$. The problem%
\begin{equation}
\left\{
\begin{array}
[c]{ll}%
-\frac{\nu_{0}}{2}\Delta_{\xi}u_{j}^{0}(\xi)+\nabla_{\xi}p_{j}^{0}(\xi)=0 &
\text{if }\xi\in\Pi^{+},\\
\mathrm{div}_{\xi}u_{j}^{0}=0 & \text{if }\xi\in\Pi^{+},\\
u_{j}^{0}=0 & \text{if }\xi_{2}=\pm\frac{1}{2},\\
u_{j}^{0}(0,\xi_{2})=\varphi_{01}\delta_{j0}-\left(
\begin{array}
[c]{c}%
u_{1,j}(0,\xi_{2})\\
u_{2,j-1}(0,\xi_{2})
\end{array}
\right)  & \text{if }\xi_{2}\in(-\frac{1}{2},\frac{1}{2}),
\end{array}
\right.  \label{B.4}%
\end{equation}
with the compatibility condition%
\begin{equation}
{\langle\varphi_{01}\delta_{j0}-u_{1,j}(0,\xi_{2})\rangle}_{2}=0, \label{B.5}%
\end{equation}
where ${\langle\psi\rangle}_{2}=\int_{-\frac{1}{2}}^{\frac{1}{2}}\psi
(x,\xi_{2})\mathrm{d}\xi_{2}$, will give the boundary layer correctors for the
velocity and for the pressure corresponding to the left end (see
\cite{Galdi}). This condition \eqref{B.5} generates a boundary condition for
$q_{j}$:%
\[
{\langle\varphi_{01}\rangle}_{2}\delta_{j0}+\frac{1}{6\nu_{0}}\left(
\frac{\partial q_{j}(0)}{\partial x_{1}}-\delta_{j0}f_{1}(0)\right)
+{\left\langle \left\{  D^{-2}\frac{\partial^{2}u_{2,j-2}}{\partial
x_{1}\partial\xi_{2}}+2D^{-2}\left(  \frac{\partial^{2}u_{1,j-2}}{\partial
x_{1}^{2}}-\frac{1}{\nu_{0}}\frac{\partial p_{j-1}}{\partial x_{1}}\right)
\right\}  _{|_{x_{1}=0}}\right\rangle }_{2}=0.
\]
In a similar way we introduce the boundary layer correctors corresponding to
the right side. The boundary layers for the velocity and pressure are defined
on $\Pi^{-}=(-\infty,0)\times(-\frac{1}{2},\frac{1}{2})$. The analogous
boundary condition for $q_{j}$ is satisfied automatically because of the
conservation of the average of $u_{1,j}(x_{1},\xi_{2})$. Indeed,%
\[
{\left\langle u_{1,j}(x_{1},\xi_{2})\right\rangle }_{2}=D^{-1}u_{1,j}%
(x_{1},\xi_{2})_{|_{\xi_{2}=\frac{1}{2}}};
\]
on the other hand, the condition
\[
{u_{2,j}}_{|_{\xi_{2}=\frac{1}{2}}}=-\frac{\partial}{\partial x_{1}}%
D^{-1}u_{1,j}(x_{1},\xi_{2})_{|_{\xi_{2}=\frac{1}{2}}}=0
\]
is equivalent to the equation on $q_{j}$, i.e. this equation on $q_{j}$ is
equivalent to the conservation law for the average ${\left\langle
u_{1,j}(x_{1},\xi_{2})\right\rangle }_{2}$. The problem satisfied by the
asymptotic solution of order $k$ is as follows:
\begin{equation}
\left\{
\begin{array}
[c]{ll}%
-\mathrm{div}(\nu(x_{1}){\mathcal{D}}\hat{u}^{k})+\nabla\hat{p}^{k}%
=f-\varepsilon^{k+1}F^{k}(x_{1},x_{2}) & \text{in }D_{\varepsilon},\\
\mathrm{div}\hat{u}^{k}=0 & \text{in }D_{\varepsilon},\\
\hat{u}_{\varepsilon}(x_{1},\frac{\varepsilon}{2})=0 & \text{for }x_{1}%
\in(0,1),\\
\hat{u}_{\varepsilon}(x_{1},-\frac{\varepsilon}{2})=0 & \text{for }x_{1}%
\in(0,1),\\
\hat{u}_{\varepsilon}(0,x_{2})=\varepsilon^{2}\varphi_{0}(\frac{x_{2}%
}{\varepsilon})+\varepsilon^{k+3}u_{2,k}(0,\xi_{2})e_{2}+\displaystyle{\sum
_{j=0}^{k}}\varepsilon^{j+2}u_{j}^{(1)}\left(  -\frac{1}{\varepsilon}%
,\frac{x_{2}}{\varepsilon}\right)  & \text{for }x_{2}\in(-\frac{\varepsilon
}{2},\frac{\varepsilon}{2}),\\
\hat{u}_{\varepsilon}(1,x_{2})=\varepsilon^{2}\varphi_{1}(\frac{x_{2}%
}{\varepsilon})+\varepsilon^{k+3}u_{2,k}(1,\xi_{2})e_{2}+\displaystyle{\sum
_{j=0}^{k}}\varepsilon^{j+2}u_{j}^{(0)}\left(  \frac{1}{\varepsilon}%
,\frac{x_{2}}{\varepsilon}\right)  & \text{for }x_{2}\in(-\frac{\varepsilon
}{2},\frac{\varepsilon}{2}),
\end{array}
\right.  \label{B.6}%
\end{equation}
Here $F^{k}$ is given by \eqref{rol}.

Mention that the boundary conditions for $u^{k}$ on $x_{1}=0,1$ are not yet
satisfied exactly because the traces of each boundary layer on the opposite
side of the rectangle  are exponentially small but do not vanish completely.
These traces should be eliminated by a small additional corrector.

Let us construct a new function $\hat{U}^{k}$ which satisfies the same
boundary conditions as $u^{k}$ on $x_{1}=0,1$.\newline Let us describe its
construction: let ${U}^{k}~:D_{\varepsilon}\longrightarrow{\mathbb{R}}^{2}$ be
a solution of the following problem:%
\begin{equation}
\left\{
\begin{array}
[c]{ll}%
U^{k}\in\left(  H^{1}(D_{\varepsilon})\right)  ^{2}, & \\
\mathrm{div}U^{k}=0 & \text{in }D_{\varepsilon},\\
{U}^{k}(x_{1},\frac{\varepsilon}{2})=0 & \text{for }x_{1}\in(0,1),\\
{U}^{k}(x_{1},-\frac{\varepsilon}{2})=0 & \text{for }x_{1}\in(0,1),\\
{U}^{k}(0,x_{2})=\varepsilon^{k+3}u_{2,k}(0,\xi_{2})e_{2}+\displaystyle{\sum
_{j=0}^{k}}\varepsilon^{j+2}u_{j}^{(1)}\left(  -\frac{1}{\varepsilon}%
,\frac{x_{2}}{\varepsilon}\right)  & \text{for }x_{2}\in(-\frac{\varepsilon
}{2},\frac{\varepsilon}{2}),\\
{U}^{k}(1,x_{2})=\varepsilon^{k+3}u_{2,k}(1,\xi_{2})e_{2}+\displaystyle{\sum
_{j=0}^{k}}\varepsilon^{j+2}u_{j}^{(0)}\left(  \frac{1}{\varepsilon}%
,\frac{x_{2}}{\varepsilon}\right)  & \text{for }x_{2}\in(-\frac{\varepsilon
}{2},\frac{\varepsilon}{2}).
\end{array}
\right.  \label{B.7}%
\end{equation}

\begin{prop}
Problem \eqref{B.7} has at least one solution with the property
\begin{equation}
\Vert U^{k}\Vert_{\left(  H^{1}(D_{\varepsilon})\right)  ^{2}}=O\left(
\varepsilon^{k+\frac{3}{2}}\right)  . \label{B.8}%
\end{equation}

\end{prop}

\textit{Proof.} Define $w^{k}_{\varepsilon}~:D \longrightarrow{\mathbb{R}}%
^{2}$, where $D=(0,1)\times(-\frac{1}{2},\frac{1}{2})$,
\[
w^{k}_{\varepsilon}(y_{1},y_{2})= \left(
\begin{array}
[c]{c}%
\varepsilon\left( U^{k}\right) _{1}\\
\left( U^{k}\right) _{2}%
\end{array}
\right) , \quad\mbox{with } (y_{1},y_{2}) = (x_{1},\frac{x_{2}}{\varepsilon
}).
\]
Obvious computations lead to the following problem for $w^{k}_{\varepsilon}$:
\begin{equation}
\left\{
\begin{array}
[c]{ll}%
\mathrm{div}_{y} w^{k}_{\varepsilon}= 0 & \quad\text{in } D,\\
w^{k}_{\varepsilon}(y_{1},\frac{1}{2})=0 & \quad\text{for } y_{1}\in(0,1),\\
w^{k}_{\varepsilon}(y_{1},-\frac{1}{2})= 0 & \quad\text{for } y_{1}\in(0,1),\\
w^{k}_{\varepsilon}(0, y_{2})= \varepsilon^{k+3} u_{2,k}(0,y_{2})e_{2}
+\left(
\begin{array}
[c]{c}%
\displaystyle{\sum_{j=0}^{k}}\varepsilon^{j+3} u^{(1)}_{1,j}\left(  -\frac
{1}{\varepsilon},y_{2}\right) \\
\displaystyle{\sum_{j=0}^{k}}\varepsilon^{j+2} u^{(1)}_{2,j}\left(  -\frac
{1}{\varepsilon},y_{2}\right)
\end{array}
\right)  & \quad\text{for } y_{2}\in(-\frac{1}{2},\frac{1}{2}),\\
w^{k}_{\varepsilon}(1, y_{2})= \varepsilon^{k+3} u_{2,k}(1,y_{2})e_{2}
+\left(
\begin{array}
[c]{c}%
\displaystyle{\sum_{j=0}^{k}} \varepsilon^{j+3} u^{(0)}_{1,j}\left(  \frac
{1}{\varepsilon},y_{2}\right) \\
\displaystyle{\sum_{j=0}^{k}}\varepsilon^{j+2} u^{(0)}_{2,j}\left(  \frac
{1}{\varepsilon},y_{2}\right)
\end{array}
\right)  & \quad\text{for } y_{2}\in(-\frac{1}{2}, \frac{1}{2}),\\
&
\end{array}
\right.
\end{equation}
Here $\mathrm{div}_{y}$ is the divergence in $y-$variables. As in
\cite{Girault-Raviart} we can prove that there exists a function
$w^{k}_{\varepsilon}\in\left( H^{1}(D_{\varepsilon})\right) ^{2}$ so that
$\|w^{k}_{\varepsilon}\|_{\left( H^{1}(D)\right) ^{2}}\leq C \|w^{k}%
_{\varepsilon}\|_{\left( H^\frac{1}{2}(\partial D)\right) ^{2}}$ with the
constant $C$ independent of $\varepsilon$. Using the properties of the
boundary layer correctors (their exponential decay rate), we get:
\[
\|w^{k}_{\varepsilon}\|_{\left( H^{1}(D)\right) ^{2}} = O\left( \varepsilon
^{k+3} \right)
\]
Direct computations give $\|U^{k}\|_{\left( H^{1}(D_{\varepsilon})\right)
^{2}} \leq\frac{1}{\varepsilon^{\frac{3}{2}}}\|w^{k}_{\varepsilon}\|_{\left(
H^{1}(D)\right) ^{2}}$. Combining these two estimates we achieve the proof.

The function%
\begin{equation}
\hat{U}^{k}=\hat{u}^{k}-U^{k} \label{B.9}%
\end{equation}
satisfies the same boundary conditions as $u$ in $x=0,1$. The problem for the
new functions $\hat{U}^{k}$, $\hat{p}^{k}$ is an obvious consequence of
\eqref{B.7} and \eqref{B.6}:%
\begin{equation}
\left\{
\begin{array}
[c]{ll}%
-\mathrm{div}(\nu(x_{1}){\mathcal{D}}\hat{U}^{k})+\nabla\hat{p}^{k}%
=f-\varepsilon^{k+1}F^{k}(x_{1},x_{2})-\mathrm{div}(\nu(x_{1}){\mathcal{D}}%
{U}^{k}) & \text{in }D_{\varepsilon},\\
\mathrm{div}\hat{U}^{k}=0 & \text{in }D_{\varepsilon},\\
\hat{U}_{\varepsilon}(x_{1},\frac{\varepsilon}{2})=0 & \text{for }x_{1}%
\in(0,1),\\
\hat{U}_{\varepsilon}(x_{1},-\frac{\varepsilon}{2})=0 & \text{for }x_{1}%
\in(0,1),\\
\hat{U}_{\varepsilon}(0,x_{2})=\varepsilon^{2}\varphi_{0}(\frac{x_{2}%
}{\varepsilon}) & \text{for }x_{2}\in(-\frac{\varepsilon}{2},\frac
{\varepsilon}{2}),\\
\hat{U}_{\varepsilon}(1,x_{2})=\varepsilon^{2}\varphi_{1}(\frac{x_{2}%
}{\varepsilon}) & \text{for }x_{2}\in(-\frac{\varepsilon}{2},\frac
{\varepsilon}{2}).
\end{array}
\right.  \label{B.10}%
\end{equation}

\begin{theorem}
Let $(\hat{u}^{k},\hat{p}^{k})$ be the asymptotic solution given by
\eqref{B.1} and $(u,p)$ the exact solution of \eqref{1.1}. Then the following
estimates hold:%
\begin{equation}
\left\{
\begin{array}
[c]{l}%
\Vert\hat{u}^{k}-u_{\varepsilon}\Vert_{{H(D_{\varepsilon})}^{2}}=O\left(
\varepsilon^{k+\frac{3}{2}}\right)  ,\\
\Vert\nabla p^{k}-\nabla p_{\varepsilon}\Vert_{H^{-1}(D_{\varepsilon}%
)}=O\left(  \varepsilon^{k+\frac{3}{2}}\right)  .
\end{array}
\right.  \label{B.11}%
\end{equation}

\end{theorem}

\textit{Proof.} From \eqref{B.9} it follows:
\[
\| \hat{u}^{k}-u_{\varepsilon}\|_{{H(D_{\varepsilon})}^{2}} \leq\|\hat{U}%
^{k}-u_{\varepsilon} \|_{{H(D_{\varepsilon})}^{2}} + \|U^{k}%
\|_{{H(D_{\varepsilon})}^{2}}=O\left( \varepsilon^{k+\frac{5}{2}} \right) +
O\left( \varepsilon^{k+\frac{3}{2}} \right) .
\]
From \eqref{rol}, \eqref{B.8} and \eqref{B.10} it follows: $\| \hat{u}%
^{k}-u\|_{{H(D_{\varepsilon})}^{2}} = O\left( \varepsilon^{k+\frac{3}{2}}
\right) $. The estimate for pressure is a consequence of $\eqref{B.11}_{1}$
and the a priori estimate.

\section{Flow in tube structures\label{sect3}}

In this section we are going to construct an asymptotic expansion to the
solution of problem \eqref{1.1}, stated in a tube structure containing one
bundle. We justify the error estimate. Let us define a tube structure
containing one bundle.

\begin{figure}
[ptb]
\begin{center}
\includegraphics[
height=1.5835in,
width=1.7832in
]%
{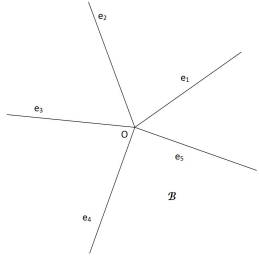}%
\caption{One bundle of segments ${\mathcal{B}}$}%
\label{f4}%
\end{center}
\end{figure}

\begin{figure}
[ptb]
\begin{center}
\includegraphics[
height=1.0559in,
width=2.7285in
]%
{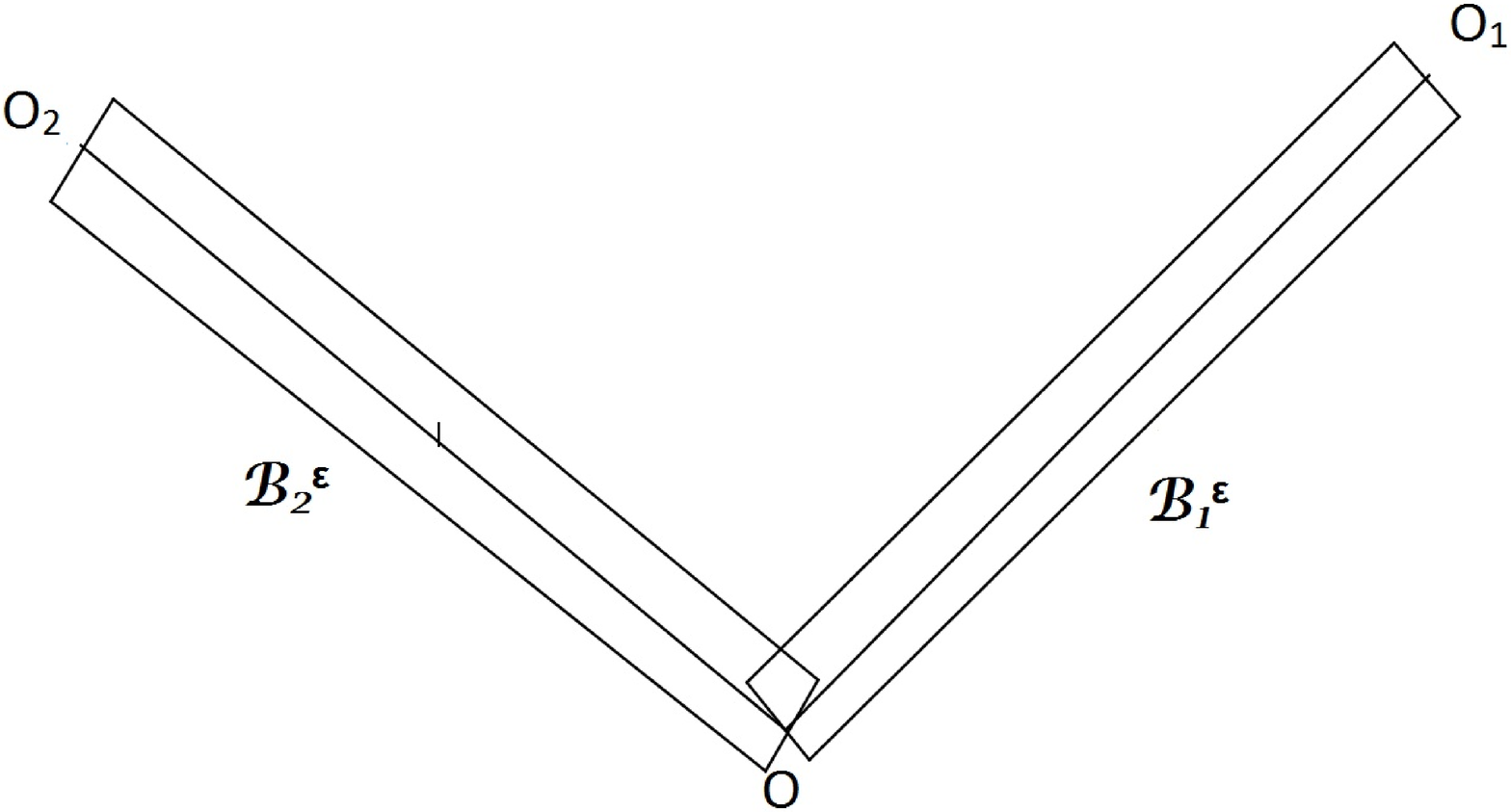}%
\caption{The rectangles ${\mathcal{B}}_{j}^{\varepsilon}$}%
\label{f5}%
\end{center}
\end{figure}

Let $e_{1},e_{2},\dots,e_{n}$ be $n$ closed segments in ${\mathbb{R}}^{2}$,
which have a single common point $O$ (i.e. the origin of the co-ordinate
system), and let it be the common end point of all these segments. Let
$\beta_{1},\beta_{2},\dots,\beta_{n}$ be $n$ bounded segments in ${\mathbb{R}%
}^{2}$ containing the point $O$, the middle point of all segments, and such
that $\beta_{j}$ is orthogonal to $e_{j}$ (for simplicity assume that the
length $\left\vert \beta_{j}\right\vert $ of each $\beta_{j}$ is equal to 1).
Let $\beta_{j}^{\varepsilon}$ be the image of $\beta_{j}$ obtained by a
homothetic contraction in $\frac{1}{\varepsilon}$ times with the center $O$.
Denote ${\mathcal{B}}_{j}^{\varepsilon}$ the open rectangles with the bases
$\beta_{j}^{\varepsilon}$ and with the heights $e_{j}$, denote also
$\hat{\beta_{j}^{\varepsilon}}$ the second base - side of each rectangle
${\mathcal{B}}_{j}^{\varepsilon}$ and let $O_{j}$ be the end of the segment
$e_{j}$ which belongs to the base $\hat{\beta_{j}^{\varepsilon}}$ (see Fig.
\ref{f5}). Define the graph of the tube structure as the bundle of segments
$e_{j}$ having a common point $O$ (see Fig. \ref{f4})%
\[
{\mathcal{B}}=\bigcup_{j=1}^{n}e_{j}.
\]
Denote below $O_{0}=O$. Let $\gamma_{j}^{\varepsilon}$, $j=0,1,\dots,n$, be
the images of the bounded domains $\gamma_{j}$ (such that $\bar{\gamma}_{j}$
contains the end of the segment $O_{j}$ and is independent of $\varepsilon$)
obtained by a homothetic contraction in $\frac{1}{\varepsilon}$ times with the
center $O_{j}$.

\begin{figure}
[ptb]
\begin{center}
\includegraphics[
height=2.2961in,
width=2.3765in
]%
{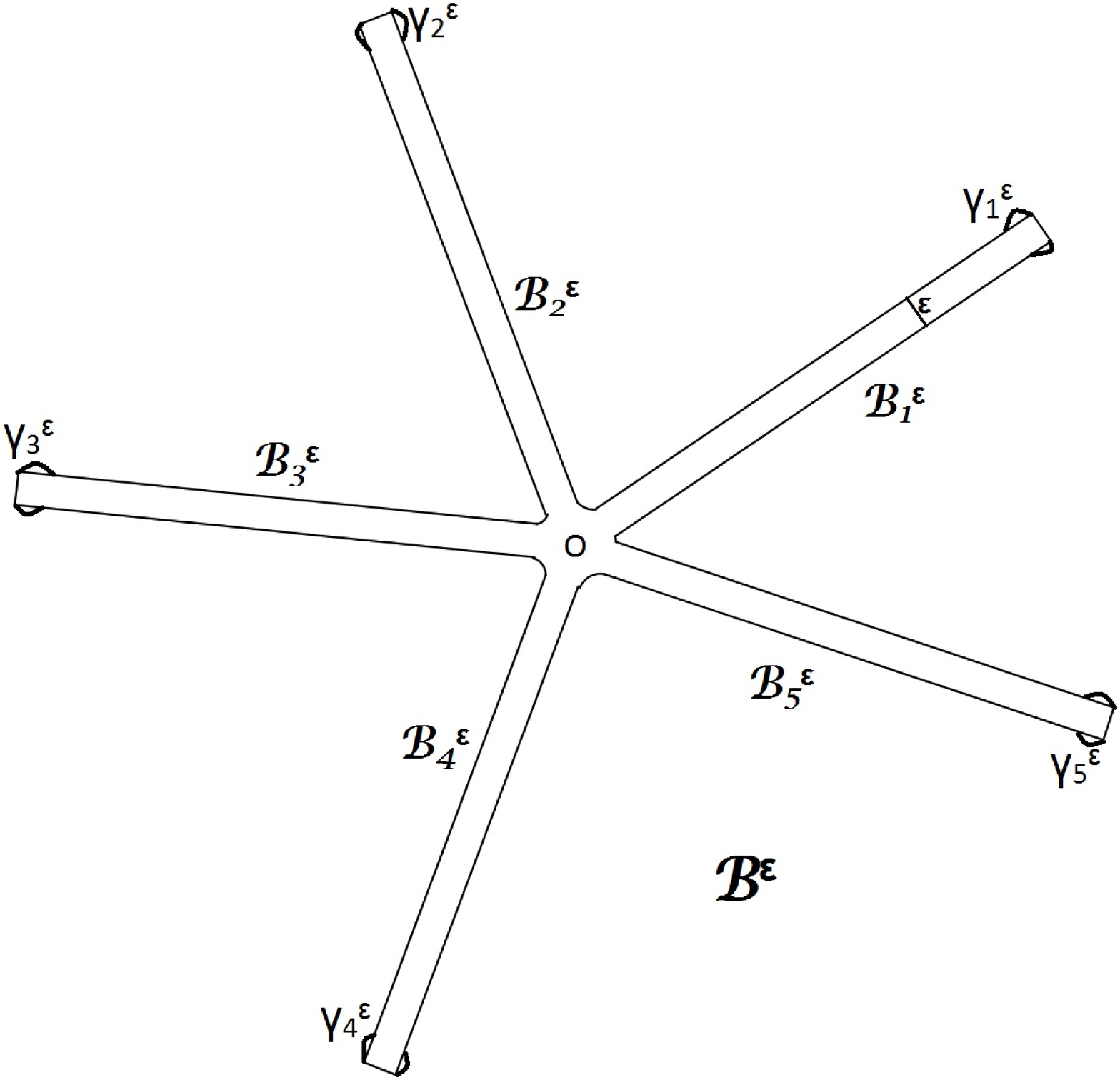}%
\caption{One bundle tubular structure ${\mathcal{B}}^{\varepsilon}$}%
\label{f3}%
\end{center}
\end{figure}

Define the tube structure associated with the bundle ${\mathcal{B}}$ as a
bounded domain (see Fig. \ref{f3}):
\[
{\mathcal{B}}^{\varepsilon}=\left( \left(  \bigcup_{j=1}^{n}\bar{\mathcal{B}%
}_{j}^{\varepsilon}\right)  \bigcup\left(  \bigcup_{j=1}^{n}\bar{\gamma}%
_{j}^{\varepsilon}\right) \right) ^{\prime}%
\]
{Here the prime stands for the set of the interior points.} Assume that
$\partial{\mathcal{B}}^{\varepsilon}\in C^{2}$ (the result may be generalized
for the case of the piecewise smooth boundary $\partial{\mathcal{B}%
}^{\varepsilon}$ with no reentrant corners). Assume that the bases $\hat
{\beta_{j}^{\varepsilon}}$ of ${\mathcal{B}}_{j}^{\varepsilon},$ $j=1,\dots
,n$, are some parts of $\partial{\mathcal{B}}^{\varepsilon}$. We add the
domains $\gamma_{j}^{\varepsilon}$, $j=0,1,\dots,n$, to smoothen the boundary
of the tube structure.%

%

%

Consider the following system of equations:%
\begin{equation}
\left\{
\begin{array}
[c]{ll}%
-\mathrm{div}(\nu(x){\mathcal{D}}u_{\varepsilon})+\nabla p_{\varepsilon
}=f(x) & \text{in }{\mathcal{B}}^{\varepsilon},\\
\mathrm{div}u_{\varepsilon}=0 & \text{in }{\mathcal{B}}^{\varepsilon},\\
u_{\varepsilon}=g & \text{on }\partial{\mathcal{B}}^{\varepsilon}.
\end{array}
\right.  \label{1}%
\end{equation}
Here, $g=0$ on the lateral boundary of the rectangles composing ${\mathcal{B}%
}^{\varepsilon}$; moreover $g=0$ anywhere with the exception of the sides
$\hat{\beta_{j}^{\varepsilon}}$ of the rectangles ${\mathcal{B}}%
_{j}^{\varepsilon}$ (these sides are assumed to belong to the boundary of the
tube structure); $g\in C^{2}(\hat{\beta}_{j}^{\varepsilon})$, and for each
$j$, $g=\varepsilon^{2}g_{j}(\frac{x-O_{j}}{\varepsilon})$ on $\hat{\beta
_{j}^{\varepsilon}}$, the vector valued functions $g_{j}\in C^{2}$ do not
depend on $\varepsilon$. Let $f$ be a vector-valued function of $(L^{2}%
({\mathcal{B}}^{\varepsilon}))^{2}$. The solvability condition gives the
relation
\begin{equation}
\int_{\partial{\mathcal{B}}^{\varepsilon}}g.n\mathrm{d}s=0. \label{2}%
\end{equation}
Introduce the local system of coordinates $Ox_{1}^{e_{j}}x_{2}^{e_{j}}$
associated with the segment $e_{j}$ such that the direction of the axis
$Ox_{1}^{e_{j}}$ coincides with the direction of the segment $OO_{j}$, i.e.
$x_{1}^{e_{j}}$ is the longitudinal coordinate. The axes $Ox_{1}^{e_{j}}%
x_{2}^{e_{j}}$ form a Cartesian coordinate system. Denote $d_{0}\varepsilon$
the infimum of radius of all circles with the center $O$ such that every point
of it belongs only to not more than one of the rectangles ${\mathcal{B}}%
_{j}^{\varepsilon}$, $j=1,\dots,n$ and $d_{1}$ is the maximal diameter of the
domains $\gamma_{0},\gamma_{1},\dots,\gamma_{n}$. We finally introduce the
notation
\[
\hat{d}_{0}\varepsilon=\max\{d_{0}\varepsilon,d_{1}\varepsilon\}.
\]
Consider the right hand side vector valued function $f$ "concentrated" in some
neighborhoods of the nodes $O_{j}$ and diffused in the rectangles, i.e.
\begin{equation}
\begin{aligned} f&=\Phi_j\left(\frac{x-O_j}{\varepsilon} \right) , &\text{for } |x-O_j| < \hat{d_0}\varepsilon, j=0,\dots,n,\\ f&=f_j(x_1^{e_j}), &\text{for } |x-O_j| > \hat{d_0}\varepsilon,\  x_1^{e_j} \in (0,|e_j|), j=1,\dots,n. \end{aligned} \label{3}%
\end{equation}
Here { $f_{j}\in C_{0}^{\infty}([0,|e_{j}|]),\Phi_{j}\in C_{0}^{1}(Q)$,}
$(j=0,1,\dots,n),$ where $Q$ is a ball $|\xi|<\hat{d_{0}}$. Assume that
\begin{equation}
\nu(x)=\nu_{0}+\nu_{j}(x_{1}^{e_{j}}) \label{4}%
\end{equation}
such that $\nu_{j}(x_{1}^{e_{j}})=0$ for all $x_{1}^{e_{j}}\in\lbrack
0,\beta]\cup\lbrack|e_{j}|-\beta;|e_{j}|]$, where $\beta$ is a positive
constant such that $\beta<min_{j}\frac{|e_{j}|}{4}$; $\nu\in C^{2}$ and there
exist $\kappa_{0}\in{\mathbb{R}}^{+}$ such that $\nu(x)>\kappa_{0}$ for all
$x\in{\mathcal{B}}^{\varepsilon}$. Without loss of generality we may assume
that $f_{j}\left(  x_{1}^{e_{j}}\right)  =0$ for all $x_{1}^{e_{j}}\in
\lbrack0,\beta]\cup\lbrack|e_{j}|-\beta;|e_{j}|]$.

Let $H_{\mathrm{div}=0}({\mathcal{B}}^{\varepsilon})$ be space of the
divergence free vector valued functions from $H^{1}({\mathcal{B}}%
^{\varepsilon})$. Let $H_{\mathrm{div}=0}^{0}({\mathcal{B}}^{\varepsilon})$ be
the subspace of vector valued functions of $H_{\mathrm{div}=0}({\mathcal{B}%
}^{\varepsilon})$ vanishing at the boundary. Assume that $g$ can be continued
in ${\mathcal{B}}^{\varepsilon}$ as a vector valued $\hat{g}$ of
$H_{\mathrm{div}=0}({\mathcal{B}}^{\varepsilon})$. The variational formulation
for \eqref{1} is as follows: find $u_{\varepsilon}\in H_{\mathrm{div}%
=0}({\mathcal{B}}^{\varepsilon})$ such that $v_{\varepsilon}=u_{\varepsilon
}-\hat{g}\in H_{\mathrm{div}=0}^{0}({\mathcal{B}}^{\varepsilon})$, and such
that it satisfies to the integral identity
\[
\int_{{\mathcal{B}}^{\varepsilon}}\nu(x){\mathcal{D}}v_{\varepsilon
}:{\mathcal{D}}{\varphi}=\int_{{\mathcal{B}}^{\varepsilon}}f\cdot\varphi
-\int_{{\mathcal{B}}^{\varepsilon}}\nu(x){\mathcal{D}}\hat{g}:{\mathcal{D}%
}\varphi,\ \ \forall\varphi\in H_{\mathrm{div}=0}^{0}({\mathcal{B}%
}^{\varepsilon}).
\]
The Riesz theorem will give the existence and the uniqueness of such a
solution because the norms $\Vert v\Vert=\sqrt{\int_{{\mathcal{B}%
}^{\varepsilon}}\nu(x){\mathcal{D}}v:{\mathcal{D}}v}$ and $\Vert
v\Vert_{(H^{1}({\mathcal{B}}^{\varepsilon}))^{2}}$ are equivalent. We have as
a consequence that
\[
\Vert v_{\varepsilon}\Vert_{(H^{1}({\mathcal{B}}^{\varepsilon}))^{2}}%
\leq{\mathcal{C}}(C_{PF},\kappa_{0})\left(  \Vert f\Vert_{L^{2}({\mathcal{B}%
}^{\varepsilon})^{2}}+\Vert\hat{g}\Vert_{(H^{1}({\mathcal{B}}^{\varepsilon
}))^{2}}\right)  ,
\]
where $C_{PF}$ is independent on $\varepsilon$ (see \cite{pana2005}) and
$\kappa_{0}$ is the lower bound of the viscosity \eqref{4}..

\begin{prop}
If $u_{\varepsilon}$ is a weak solution for problem \eqref{1} then there
exists a distribution $p_{\varepsilon}\in{{\mathcal{D}}}^{^{\prime}%
}({\mathcal{B}}^{\varepsilon})$ such that $(u_{\varepsilon},p_{\varepsilon})$
satisfies $\eqref{1}_{1}$ in the sense of distributions. The following
inequality holds in the case of $\hat{g}=0$:
\[
\Vert\nabla p_{\varepsilon}\Vert_{H^{-1}(B^{\varepsilon})}\leq C\Vert
f\Vert_{\left(  L^{2}(B^{\varepsilon})\right)  ^{2}}.
\]
Here $C$ is a constant independent of $\varepsilon.$
\end{prop}

The proof is similar to that of the Propositions of section \ref{sect2}.

\subsection{Asymptotic expansion\label{sect3.1}}

We construct the main part of the asymptotic expansion in a form%
\begin{equation}
u^{a}=\displaystyle{\sum_{l=0}^{k}\varepsilon^{l+2}\left\{  \sum
_{e=e_{j};j=1,\dots,n}u_{l}^{e}(x^{e,L})\chi_{\varepsilon}(x)+\sum_{i=0}%
^{n}u_{l}^{BLO_{i}}\left(  \frac{x-O_{i}}{\varepsilon}\right)  \right\}  }
\label{5}%
\end{equation}%
\begin{equation}%
\begin{array}
[c]{lcl}%
p^{a} & = & \displaystyle{\sum_{l=0}^{k}\varepsilon^{l+1}\left\{
\sum_{e=e_{j};j=1,\dots,n}p_{l}^{e}(x^{e,L})\chi_{\varepsilon}(x)+\sum
_{i=0}^{n}p_{l}^{BLO_{i}}\left(  \frac{x-O_{i}}{\varepsilon}\right)  \right\}
}+\\
&  & +\displaystyle{\sum_{e=e_{j};j=1,\dots,n}\sum_{l=0}^{k}\varepsilon
^{l}q_{l}^{e}(x_{1}^{e})\chi_{\varepsilon}(x)+\sum_{i=0}^{n}q_{0}^{e_{i}%
}(x_{1}^{e_{i}}=0)\left(  1-\chi_{\varepsilon}(x)\right)  \theta_{i}(x)}%
\end{array}
\label{6}%
\end{equation}
$x^{e,L}=(x_{1}^{e},\frac{x_{2}^{e}}{\varepsilon})$. Later, in the end of the
section we will add an exponentially small corrector multiplying the boundary
layers by a cut-off function $\eta$ in the subdomain where the boundary layers
are just exponentially small (see (\ref{15}), (\ref{16})). The last sum in
(\ref{1}) is taken for all nodes $O_{i}$ and the value $q_{0}^{e_{i}}%
(x_{1}^{e_{i}}=0)$ is calculated at the point $x=O_{i}$; function
$q_{0}^{e_{i}}$ is supposed to be continuous on the graph ${\mathcal{B}}$.
Here $\chi_{\varepsilon}(x)$ is a function equal to zero at the distance less
than $(\hat{d_{0}}+1)\varepsilon$ from $O_{j}$, $j=0,1,\dots,n$, equal  to
zero on the rectangle ${\mathcal{B}}_{j}^{\varepsilon}$ if $x_{1}^{e_{j}}%
\leq(\hat{d_{0}}+1)\varepsilon$ or if $|x_{1}^{e_{j}}-|e_{j}||\leq(\hat{d_{0}%
}+1)\varepsilon$; we suppose that function $\chi_{\varepsilon}$ is equal to
one on this rectangle if $x_{1}^{e_{j}}\geq(\hat{d_{0}}+2)\varepsilon$ and
$|x_{1}^{e_{j}}-|e_{j}||\geq(\hat{d_{0}}+2)\varepsilon$, and we define
$\chi_{\varepsilon}$ by the relations {$\chi_{\varepsilon}(x)=\chi\left(
\frac{x_{1}^{e_{j}} }{\varepsilon}\right)  $ if $(\hat{d_{0}}+1)\varepsilon
\leq x_{1}^{e_{j}}\leq(\hat{d_{0}}+2)\varepsilon$ and $\chi_{\varepsilon
}(x)=\chi\left( \frac{x_{1}^{e_{j}}-|e_{j}| }{\varepsilon}\right)  $ if
$(\hat{d_{0}}+1)\varepsilon\leq|e_{j}|-x_{1}^{e_{j}}\leq(\hat{d_{0}%
}+2)\varepsilon$}. Here $\chi$ is a differentiable on ${\mathbb{R}}$ function
of one variable, it is independent of $\varepsilon$, it is equal to zero on
the segment $[-(\hat{d_{0}}+1);(\hat{d_{0}}+1)]$ and it is equal to one on the
union of intervals $(-\infty,-(\hat{d_{0}}+2))\cup((\hat{d_{0}}+2),+\infty)$.
Moreover $\chi_{\varepsilon}$ is equal to zero on every $\gamma_{j}%
^{\varepsilon}$. The functions ${\theta_{i}}$, $i=1,...,n,$ are defined as
follows%
\[
\begin{aligned}
\theta_j=0 &\quad |x-O_j|>min_i\frac{|e_i|}{2}\\
\theta_j=1 &\quad |x-O_j|\leq min_i\frac{|e_i|}{2}\\
\end{aligned}
\]

The relation between the vector-columns $x^{T}$ and $x^{e_{j},T}$  (here $T$
is the transposition symbol) is given by
\[
x^{T}=\Gamma_{j}x^{e_{j},T}+O\qquad j=1,\dots,n
\]
where $\Gamma_{j}$ is an orthogonal matrix of passage from the canonic base to
the local one. Then applying the results of section \ref{sect2}, for every
channel ${\mathcal{B}}_{j}^{\varepsilon},$ we get $u_{l}^{e}$, $p_{l}^{e}$ and
$q_{l}^{e}$ defined up to the scalar constants $c_{l}^{e}$, $d_{l}^{e}$.
Indeed, denote $\hat{q}_{l}^{e_{j}}(x_{1}^{e_{j}})$ the solution of equation
(\ref{A.5}) with $\dfrac{\partial\hat{q}_{l}^{e_{j}}}{\partial x_{1}^{e_{j}}%
}(0)=0.$ Then the general solution of equation (\ref{A.5}) has a form:
\begin{equation}
q_{l}^{e_{j}}(x_{1}^{e_{j}})=\hat{q}_{l}^{e_{j}}(x_{1}^{e_{j}})+c_{l}^{e_{j}%
}\int_{0}^{x_{1}^{e_{j}}}\nu(s)\mathrm{d}s+d_{l}^{e_{j}}, \label{60}%
\end{equation}
where $c_{l}^{e_{j}}$ and $d_{l}^{e_{j}}$ are the undetermined constants;
\[
u_{l}^{e_{j}}(x^{e_{j},L})=\Gamma_{j}\left(  \tilde{u}_{l}^{e_{j}}(x^{e_{j}%
,L})\right)  ^{T},
\]
where the second component of $\tilde{u}_{l}^{e_{j}}$ does not depend on
$c_{l}^{e_{j}}$, $d_{l}^{e_{j}}$ (see ($\ref{A.4}$)$_{2}$); the same property
holds for $p_{l}^{e_{j}}$ (see ($\ref{A.4}$)$_{3}$); the first component
$\tilde{u}_{1,l}^{e_{j}}$ (see ($\ref{A.4}$)$_{1}$) depends on $c_{l}^{e_{j}}%
$:%
\begin{align*}
\tilde{u}_{1,l}^{e_{j}}  &  =-D^{-2}\left\{  \frac{\partial^{2}\tilde
{u}_{2,l-2}^{e_{j}}}{\partial\xi_{2}^{e_{j}}\partial x_{1}^{e_{j}}}+\frac
{2}{\nu(x_{1}^{e_{j}})}\left(  \frac{\partial}{\partial x_{1}^{e_{j}}}\left(
\nu(x_{1}^{e_{j}})\frac{\partial\tilde{u}_{1,l-2}^{e_{j}}}{\partial
x_{1}^{e_{j}}}\right)  -\frac{\partial p_{l-1}^{e_{j}}}{\partial x_{1}^{e_{j}%
}}\right)  \right\}  +\\
&  +\frac{2}{\nu(x_{1}^{e_{j}})}N_{1}(\xi_{2}^{e_{j}})\left(  \frac
{\partial\hat{q}_{l}^{e_{j}}}{\partial x_{1}^{e_{j}}}-f_{1}\delta_{j0}\right)
+2c_{l}^{e_{j}}N_{1}(\xi_{2}^{e_{j}})\\
&  =\hat{u}_{1,l}^{e_{j}}+2c_{l}^{e_{j}}N_{1}(\xi_{2}^{e_{j}}).
\end{align*}
For $x_{1}^{e_{j}}\in\lbrack0,{\beta}]\bigcup[|e_{j}|-{\beta},|e_{j}|]$,
$\hat{u}_{l}^{e_{j}}$, $\tilde{p}_{l}^{e_{j}}$ and $\partial\hat{q}_{l}%
^{e_{j}}/\partial x_{1}^{e_{j}}$ may be taken equal to zero because the right
hand side $f_{j}$ is equal to zero for this values of $x_{1}^{e_{j}}$, while
the flow rate $\int_{\beta_{j}}\tilde{u}_{1,l}^{e_{j}}(\xi_{2}^{e_{j}}%
)d\xi_{2}^{e_{j}}$ is constant on $e_{j}$.{ $N_{1}$ is the function introduced
in section 2.
}

To get the problems for the boundary layers, we introduce the domain
$\Omega_{O_{0}}=\displaystyle{\cup_{j=1}^{n}\tilde{\Omega}_{j}\cup\gamma_{0}}%
$, where $\tilde{\Omega}_{j}$ are the half-infinite strips obtained from
${\mathcal{B}}_{j}^{\varepsilon}$ by infinite extension behind the base
$\tilde{\beta}_{j}^{\varepsilon}$ and by homothetic dilatation in $\frac
{1}{\varepsilon}$ times (with respect to the point $O$); let $\Omega_{j}$ be
obtained from $\tilde{\Omega}_{j}$ by a symmetric reflection relatively to the
line containing $\beta_{j}^{\varepsilon}$ and let $\Omega_{O_{j}}%
=\tilde{\Omega}_{j}\cup\gamma_{j}^{t}$, where $\gamma_{j}^{t}$ is obtained
from $\gamma_{j}$ by a translation (such that the point $O_{j}$ becomes $O$).

Since $\nu(x)=\nu_{0}$ for all $x_{1}^{e_{j}}\in\lbrack0,\beta]\cup
\lbrack|e_{j}|-\beta;|e_{j}|]$, the boundary layer solution is a pair
constituted of a vector valued function $u_{l}^{BLO_{j}}$ and a scalar
function $p_{l}^{BLO_{j}}$ satisfying to the Stokes system:%
\begin{equation}
\left\{
\begin{array}
[c]{l}%
-\frac{\nu_{0}}{2}\Delta_{\xi}u_{l}^{BLO_{0}}+\nabla_{\xi}p_{l}^{BLO_{0}}%
=\Phi_{0}(\xi)\delta_{l,0}+\\
\ \ \ +\displaystyle{\sum_{e=e_{j};j=1,\dots,n}}\left\{  c_{l}^{e}\left(
\nu_{0} \Delta_{\xi}\left(  \chi_{j}(\xi_{1}^{e})\Gamma_{j}\left(  N_{1}%
(\xi_{2}^{e}),0\right)  ^{T}\right)  -\nabla_{\xi}\left(  \chi_{j}(\xi_{1}%
^{e})\xi_{1}^{e}\right)  \right)  -d_{l+1}^{e}\nabla_{\xi}\chi_{j}(\xi_{1}%
^{e})\right\}  ,\\
\\
\mathrm{div}_{\xi}u_{l}^{BLO_{0}}=-\displaystyle{\sum_{e=e_{j};\ j=1,\dots,n}%
}\mathrm{div}_{\xi}\left(  c_{l}^{e}\left(  \chi_{j}(\xi_{1}^{e})\Gamma
_{j}\left(  N_{1}(\xi_{2}^{e}),0\right)  ^{T}\right)  \right)  ,\text{ \ \ if
}\xi\in\Omega_{O_{0}},\\
{u_{l}^{BLO_{0}}}{|_{\partial\Omega_{O_{0}}}}=0,
\end{array}
\right.  \label{8}%
\end{equation}
and for $j=1,\dots,n,$%
\begin{equation}
\left\{
\begin{array}
[c]{l}%
-\frac{\nu_{0}}{2}\Delta_{\hat{\xi}}u_{l}^{BLO_{j}}+\nabla_{\hat{\xi}}%
p_{l}^{BLO_{j}}=\\
\ \ \ =\Phi_{j}(\hat{\xi})\delta_{l,0}+\hat{c}_{l}^{e}\left(  \nu_{0}
\Delta_{\hat{\xi}}\left(  \chi_{j}(\hat{\xi}_{1}^{e})\right)  \hat{\Gamma}%
_{j}\left(  N_{1}(\hat{\xi}_{2}^{e}),0\right)  ^{T}-\nabla_{\hat{\xi}}\left(
\chi_{j}(\hat{\xi}_{1}^{e})\hat{\xi}_{1}^{e}\right)  \right)  -\hat{d}%
_{l+1}^{e}\nabla_{\hat{\xi}}\chi_{j}(\hat{\xi}_{1}^{e}),\\
\\
\mathrm{div}_{\hat{\xi}}u_{l}^{BLO_{j}}=-\hat{c}_{l}^{e_{j}}\mathrm{div}%
_{\hat{\xi}}\left(  \chi_{j}(\hat{\xi}_{1}^{e_{j}})\hat{\Gamma}_{j}\left(
N_{1}(\hat{\xi}_{2}),0\right)  ^{T}\right)  ,\text{ \ \ if }\hat{\xi}\in
\Omega_{O_{j}},\\
{u_{l}^{BLO_{j}}}_{|_{\partial\Omega_{O_{j}},\hat{\xi}_{1}^{e_{j}}=0}}%
=g_{j}\delta_{l,0},\\
{u_{l}^{BLO_{j}}}_{|_{\partial\Omega_{O_{j}},\hat{\xi}_{1}^{e_{j}}\neq0}}=0.
\end{array}
\right.  \label{9}%
\end{equation}
The variable $\hat{\xi}_{1}^{e_{j}}$ is opposite to $\xi_{1}^{e_{j}}$, i.e. to
the first component of the vector $\Gamma_{j}^{T}\xi^{T}$. So $\hat{\xi}%
_{1}^{e_{j}}=\hat{\Gamma}_{j}^{T} \xi^{T}$, where $\hat{\Gamma}_{j}=\hat
{I}d\Gamma_{j}$ and $\hat{I}d$ is the diagonal matrix with the diagonal
elements $-1,1$. The constants $\hat{c}_{l}^{e_{j}}$, $\hat{d}_{l}^{e_{j}}$
are defined in such a way that the functions $c_{l}^{e_{j}}\int_{0}%
^{x_{1}^{e_{j}}}\nu(s)\mathrm{d}s+d_{l}^{e_{j}}$ and $\hat{c}_{l}^{e_{j}}%
\int_{0}^{|e_{j}|-x_{1}^{e_{j}}}\nu(s)\mathrm{d}s+\hat{d}_{l}^{e_{j}}$ are
equal, i.e.%
\begin{equation}
c_{l}^{e_{j}}=-\hat{c}_{l}^{e_{j}},\ \ \ \hat{d}_{l}^{e_{j}}=c_{l}^{e_{j}}%
\int_{0}^{|e_{j}|}\nu(s)\mathrm{d}s+d_{l}^{e_{j}}. \label{10}%
\end{equation}
Assume that every term in the sum $\displaystyle{\sum_{e=e_{j};j=1,\dots,n}}$
in \eqref{8} is defined only in the branch of $\Omega_{O_{0}}$, corresponding
to $e=e_{j}$, and it vanishes in $\gamma_{0}$.

The solutions of these boundary layer problems decay exponentially  at
infinity and the constants $c_{l}^{e_{j}}$, $\hat{c}_{l}^{e_{j}}$, ${d}%
_{l}^{e_{j}}$ and $\hat{d}_{l}^{e_{j}}$ are chosen from the conditions of
existence of such solutions (see \cite{Naza}). Let us define first $\hat
{c}_{l}^{e_{j}}$ from the condition of exponential decaying of $u_{l}%
^{BLO_{j}}$ at infinity:%
\[
\int_{\Omega_{O_{j}}}\hat{c}_{l}^{e_{j}}\mathrm{div}_{\hat{\xi}}\left(
\chi_{j}(\hat{\xi}_{1}^{e_{j}})\hat{\Gamma}_{j}\left(  N_{1} ({\xi}_{2}%
^{e}),0\right)  ^{T}\right)  d\hat{\xi}=\int_{\beta_{j}}\left(  \hat{\Gamma
}_{j}^{T}g_{j}\right)  ^{1}\mathrm{d}{\xi}_{2}^{e}\delta_{l,0},
\]
i.e.
\begin{equation}
-\int_{\beta_{j}}N_{1}(\xi_{2}^{e_{j}})\mathrm{d}{\xi}_{2}^{e}\hat{c}%
_{l}^{e_{j}}=\int_{\beta_{j}}\left(  \hat{\Gamma}_{j}^{T}g_{j}\right)
^{1}\mathrm{d}{\xi}_{2}^{e}\delta_{l,0}, \label{11}%
\end{equation}
where the upper index $1$ corresponds to the first component of the vector.

Then we find $\hat{c}_{l}^{e_{j}}$, and $\hat{d}_{l}^{e_{j}}$ as defined in
\eqref{10}. Then we determine the constants ${d}_{l+1}^{e_{j}}$ from the
condition of the exponential decaying of $p_{l}^{BLO_{0}}$ at infinity. To
this end, consider first problem \eqref{8} without the last term in equation
${\eqref{8}}_{1}$, i.e.%
\begin{equation}
\left\{
\begin{array}
[c]{l}%
-\frac{\nu_{0}}{2}\Delta_{\xi}\bar{u}_{l}^{BLO_{0}}+\nabla_{\xi}\bar{p}%
_{l}^{BLO_{0}}=\\
\ \ \ \ =\Phi_{0}(\xi)\delta_{l,0}+\displaystyle{\sum_{e=e_{j};j=1,\dots,n}%
}\left\{  c_{l}^{e}\left(  \nu_{0}\Delta_{\xi}\left(  \chi_{j}(\xi_{1}%
^{e})\right)  \Gamma_{j}\left(  N_{1}(\xi_{2}^{e}),0\right)  ^{T}-\nabla_{\xi
}\left(  \chi_{j}(\xi_{1}^{e})\xi_{1}^{e}\right)  \right)  \right\}  ,\\
\\
\mathrm{div}_{\xi}\bar{u}_{l}^{BLO_{0}}=-\displaystyle{\sum_{e=e_{j}%
;j=1,\dots,n}}\mathrm{div}_{\xi}c_{l}^{e}\left(  \chi_{j}(\xi_{1}^{e}%
)\Gamma_{j}\left(  N_{1}(\xi_{2}^{e}),0\right)  ^{T}\right)  ,\quad\xi
\in\Omega_{O_{0}},\\
\\
{\bar{u}_{l}^{BLO_{0}}}{|_{\partial{\Omega}_{O_{0}}}}=0.
\end{array}
\right.  \label{12}%
\end{equation}
Here the constants $c_{l}^{e_{j}}$ are just defined by \eqref{10} and
\eqref{11} and satisfy the condition
\[
\int_{\Omega_{O_{0}}}{\sum_{e=e_{j};\ j=1,\dots,n}}\mathrm{div}_{\xi}c_{l}%
^{e}\left(  \chi_{j}(\xi_{1}^{e})\Gamma_{j}\left(  N_{1}({\xi}_{2}%
^{e}),0\right)  ^{T}\right)  d{\xi}_{2}^{e}=0
\]
i.e.
\begin{equation}
{\sum_{e=e_{j};j=0,1,\dots,n}}\int_{\beta_{j}}c_{l}^{e_{j}} N_{1}(\xi_{2}%
^{e})d{\xi}_{2}^{e}=0. \label{13}%
\end{equation}
Indeed, the choice of constants $c_{l}^{e_{j}}=-\hat{c}_{l}^{e_{j}}$ and
$c_{l}^{e_{j}}$ from \eqref{11} and condition \eqref{2} give relation \eqref{13}.

It is known that there exists the unique solution $\{\bar{u}_{l}^{BLO_{0}%
},\bar{p}_{l}^{BLO_{0}}\}$ of this problem such that $\bar{u}_{l}^{BLO_{0}}$
stabilizes to zero at infinity on every branch of $\Omega_{0}$ and $\bar
{p}_{l}^{BLO_{0}}$ stabilizes on every branch of $\Omega_{0}$ associated to
$e_{j}$, to its own constant $\bar{p}_{l}^{BLO_{0}\infty j}$. These constants
are defined uniquely up to one common additive constant, which we fix here by
a condition $\bar{p}_{l}^{BLO_{0}\infty1}=0$. Then we define%
\begin{equation}
\left\{
\begin{array}
[c]{l}%
{d}_{l+1}^{e_{j}}=-\bar{p}_{l}^{BLO_{0}\infty j},\\
\\
u_{l}^{BLO_{0}}=\bar{u}_{l}^{BLO_{0}},\\
\\
{p}_{l}^{BLO_{0}}=\bar{p}_{l}^{BLO_{0}}+\displaystyle{\sum_{e=e_{j}%
;j=1,\dots,n}}{d}_{l+1}^{e_{j}}\chi_{j}(\xi_{1}^{e_{j}})
\end{array}
\right.  \label{14}%
\end{equation}
on every branch of $\Omega_{0}$, associated with $e_{j}$, i.e. ${p}%
_{l}^{BLO_{0}}=\bar{p}_{l}^{BLO_{0}}-{\sum_{e=e_{j};j=1,\dots,n}}\bar{p}%
_{l}^{BLO_{0}\infty j}\chi_{j}.$

Obviously, this pair $\{{u}_{l}^{BLO_{0}},{p}_{l}^{BLO_{0}}\}$ satisfies
\eqref{8}. The boundary layer functions ${u}_{l}^{BLO_{j}}$ and ${p}%
_{l}^{BLO_{j}}$; $j=0,1,\dots,n$ are not defined in the vicinity of $O$.
Therefore we should change a little bit the formulas of $u^{a}$ and $p^{a}$
far from the nodes $O_{j}$, $j=0,1,\dots,n$.

Let $\eta_{j}(x_{1}^{e_{j}})$ be a smooth function defined on each segment
$e_{j}$, let it be one if $\left\vert x_{1}^{e_{j}}-\frac{|e_{j}|}%
{2}\right\vert \geq\frac{|e_{j}|}{2}-\beta$ and let it be zero if $\left\vert
x_{1}^{e_{j}}-\frac{|e_{j}|}{2}\right\vert \leq\frac{|e_{j}|}{8}$. Let
$\eta(x)=\eta_{j}(x_{1}^{e_{j}})$ for each rectangle ${\mathcal{B}}%
_{j}^{\varepsilon}$ and let $\eta=1$ on each $\gamma_{j}^{\varepsilon}$. Set
$\eta^{(j)}(x)=\eta(x)$ on $\gamma_{j}^{\varepsilon}$ and all half-rectangles
having common points with $\gamma_{j}^{\varepsilon}$ and extend it by zero on
the remaining part of ${\mathcal{B}}_{j}^{\varepsilon}$. Then we define
$u^{a}$ and $p^{a}$ as
\begin{equation}
\bar{u}^{a}=\displaystyle{\sum_{l=0}^{k}\varepsilon^{l+2}\left\{
\sum_{e=e_{j};j=0,1,\dots,n}u_{l}^{e}(x^{e,L})\chi_{\varepsilon}(x)+\sum
_{i=0}^{n}u_{l}^{BLO_{i}}\left(  \frac{x-O_{i}}{\varepsilon}\right)
\eta^{(i)}(x)\right\}  ,} \label{15}%
\end{equation}%
\begin{equation}%
\begin{array}
[c]{lcl}%
\bar{p}^{a} & = & \displaystyle{\sum_{l=0}^{k}\varepsilon^{l+1}\left\{
\sum_{e=e_{j};j=0,1,\dots,n}p_{l}^{e}(x^{e,L})\chi_{\varepsilon}(x)+\sum
_{i=0}^{n}p_{l}^{BLO_{i}}\left(  \frac{x-O_{i}}{\varepsilon}\right)
\eta^{(i)}(x)\right\}  }\\
&  & +\displaystyle{\sum_{l=0}^{k}\varepsilon^{l}q_{l}^{e}(x_{1}^{e}%
)\chi_{\varepsilon}(x)+\sum_{i=0}^{n}q_{0}^{e_{i}}(x_{1}^{e_{i}}=0)\left(
1-\chi_{\varepsilon}(x)\right)  \theta_{i}(x).}%
\end{array}
\label{16}%
\end{equation}
Let us mention that the last term (sum) in (\ref{16}) corresponds to the
boundary layer function $p_{l}^{BLO_{j}}$ for $l=-1.$

\subsection{Error estimate\label{sect3.2}}

In this section we estimate the error between the exact solution and the
asymptotic one. {Substituting the asymptotic expansions \eqref{15},
\eqref{16}, into \eqref{1}, we get the relations
\begin{equation}
\left\{
\begin{array}
[c]{ll}%
-\mathrm{div}(\nu(x){\mathcal{D}}\bar{u}^{a})+\nabla\bar{p}^{a}=f(x) & \\
\ \ \ \ \ \ -\displaystyle{\sum_{e=e_{j};j=1,\dots,n}}\left\{  {\varepsilon
}^{k+1}\Gamma_{j}\left(  F_{e_{j}}^{k}\right)  ^{T}\chi_{\varepsilon
}(x)-\varepsilon^{k}(\nabla_{\xi}\chi_{j}(\xi_{1}^{e_{j}})d_{k+1}^{e_{j}%
}\right.  & \\
\qquad\qquad\left. -\nabla_{\tilde{\xi}}\chi_{j}(\tilde{\xi}_{1}^{e_{j}}%
){d}_{k+1}^{e_{j}})\right\}  +\Psi, & \text{in }{\mathcal{B}}^{\varepsilon},\\
\mathrm{div}\bar{u}^{a}=\psi & \text{in }{\mathcal{B}}^{\varepsilon},\\
\bar{u}^{a}=g & \text{on }\partial{\mathcal{B}}^{\varepsilon},
\end{array}
\right.  \label{17}%
\end{equation}
where $F_{e_{j}}^{k}$ is the residual described in \eqref{rol}, and $\psi$ is
defined by
\begin{equation}
\label{18}\psi(x,t)=\left\{
\begin{array}
[c]{ll}%
\displaystyle 0 & \mathrm{in}~{\mathcal{B}}^{\varepsilon}\cap\left\{ x^{e_{i}%
}_{1}<\beta\right\} ,\vspace{0.1cm}\\
\displaystyle-\nabla_{x}\eta_{i}(x_{1}^{e_{i}}). u^{(BLO_{0})}\left(
\frac{x-O_{0}}{\varepsilon}\right)  & \mathrm{in}~{\mathcal{B}}^{\varepsilon
}\cap\left\{ \beta<x_{1}^{e_{i}}<\frac{3|e_{i}|}{8}\right\} ,\vspace{0.1cm}\\
0 & \mathrm{in}~{\mathcal{B}}^{\varepsilon}\cap\left\{ \frac{3|e_{i}|}%
{8}<x_{1}^{e_{i}}<\frac{5|e_{i}|}{8}\right\} ,\vspace{0.1cm}\\
\displaystyle-\nabla_{x}\eta_{i}(x_{1}^{e_{i}}).u^{(BLO_{i})}\left(
\frac{x-O_{i}}{\varepsilon}\right)  & \mathrm{in}~{\mathcal{B}}^{\varepsilon
}\cap\left\{ \frac{5|e_{i}|}{8}<x_{1}^{e_{i}}<|e_{i}|\right\} ,\vspace
{0.1cm}\\
0 & \mathrm{in}~\gamma_{i}^{\varepsilon},i=1,\dots,n,\vspace{0.1cm}\\
&
\end{array}
\right.
\end{equation}
$\Vert\Psi\Vert_{\left( L^{2}(({\mathcal{B}}^{\varepsilon})\right)  ^{2}%
}=O\left(  \exp(\frac{-c}{\varepsilon})\right) $, $\Vert\psi\Vert
_{H^{1}({\mathcal{B}}^{\varepsilon})}=O\left(  \exp(\frac{-c}{\varepsilon
})\right)  $ with a positive constant $c$ and ${\displaystyle\int
_{{\mathcal{B}}^{\varepsilon}}} \psi\mathrm{d}s=0,$ because
$\displaystyle{\int_{\partial{\mathcal{B}}^{\varepsilon}}(\bar{u}%
^{a},n)\mathrm{d}s=\int_{\partial{\mathcal{B}}^{\varepsilon}}(g,n)\mathrm{d}%
s=0}$. The exponentially decaying residuals $\Psi$ and $\psi$ appear from the
truncation of the boundary layer terms by the function $\eta$: it is different
from 1 in the part of the domain ${\mathcal{B}}^{\varepsilon}$ where
$u_{l}^{BLO_{i}}\left(  \frac{x-O_{i}}{\varepsilon}\right)  $, $p_{l}%
^{BLO_{i}}\left(  \frac{x-O_{i}}{\varepsilon}\right)  $ and their derivatives
are exponentially small. We are going to prove the estimate
\[
\Vert u-\bar{u}^{a}\Vert_{H^{1}({\mathcal{B}}^{\varepsilon})}=O\left(
{{\varepsilon}^{k+\frac{1}{2}}}\right)  .
\]
We can not apply directly the a priori estimates because $\bar{u}^{a}$ is not
divergence free.
\newline Let us construct a function $\hat{U}^{a}: {\mathcal{B}}^{\varepsilon
}\rightarrow{\mathbb{R}}^{2}$ satisfying the following properties:
\begin{equation}
\left\{
\begin{array}
[c]{ll}%
\displaystyle \hat{U}^{a}\in\left( H_{0}^{1}({\mathcal{B}}^{\varepsilon
})\right) ^{2}, & \vspace{0.1cm}\\
\displaystyle \mathrm{div}_{x} \hat{U}^{a}=\psi & \mathrm{in}\qquad
{\mathcal{B}}^{\varepsilon},\vspace{0.1cm}\\
\hat{U}^{a}=0 & \mathrm{in}\quad{\mathcal{B}}^{\varepsilon}\cap\{x^{e}%
_{1}<\hat{d}_{0}\varepsilon\}\bigcup\gamma_{j}^{\varepsilon},j=1,\dots
,n\vspace{0.1cm}\\
&
\end{array}
\right. \label{19}%
\end{equation}
}

\begin{prop}
Problem \eqref{19} has at least one solution, satisfying
\[
\Vert\hat{U}^{a} \Vert_{(H^{1}({\mathcal{B}}^{\varepsilon}))^{2}%
}=O(exp(-c/\varepsilon))
\]

\end{prop}

\textit{Proof.} Due to \eqref{19}$_{3}$ we can consider the problem \eqref{19}
as separate problem on each ${\mathcal{B}}^{\varepsilon}_{j}$. Denote by
$\hat{U}^{e_{j}}$ the restriction of $\hat{U}^{a}$ on ${\mathcal{B}%
}^{\varepsilon}_{j}$ , obviously $\hat{U}^{e_{j}}(x_{1}^{e_{j}},x_{2}^{e_{j}%
})=0$ for $(x \in{\mathcal{B}}^{\varepsilon}_{j} $ such that $x_{1}^{e_{j}%
}<\hat{d}_{0}\varepsilon$. For all $(x\in{\mathcal{B}}^{\varepsilon}_{j} $
such that $x_{1}^{e_{j}}>\hat{d}_{0}\varepsilon$ introduce the new variable
$(y_{1}^{e_{j}},y_{2}^{e_{j}})=\left( \frac{x_{1}^{e_{j}}-\hat{d}%
_{0}\varepsilon}{|e_{j}|-\hat{d}_{0}\varepsilon},\frac{x_{2}^{e_{j}}%
}{\varepsilon}\right)  $; obviously $(y_{1}^{e_{j}},y_{2}^{e_{j}}%
)\in(0,1)\times(-\frac{1}{2},\frac{1}{2}) $. Define a new function
$\mu_{\varepsilon}^{e_{j}}:(0,1)\times(-\frac{1}{2},\frac{1}{2}%
)\longrightarrow{\mathbb{R}}^{2}$, by $\mu_{\varepsilon}^{e_{j}}(y_{1}^{e_{j}%
},y_{2}^{e_{j}})= \left( \frac{1}{|e_{j}|-\hat{d}_{0}\varepsilon}\hat{U}%
_{1}^{e_{j}}((|e_{j}|-\hat{d}_{0}\varepsilon)y_{1}^{e_{j}}+ \hat{d}%
_{0}\varepsilon,\varepsilon y_{2}^{e_{j}} ),\frac{1}{\varepsilon}\hat{U}%
_{2}^{e_{j}}((|e_{j}|-\hat{d}_{0}\varepsilon)y_{1}^{e_{j}}+ \hat{d}%
_{0}\varepsilon,\varepsilon y_{2}^{e_{j}} )\right)  $. Obvious computations
lead to the following problem for $\mu_{\varepsilon}^{e_{j}}$:
\begin{equation}
\left\{
\begin{array}
[c]{ll}%
\displaystyle \mathrm{div}_{y} \mu_{\varepsilon}^{e_{j}}=\psi((|e_{j}|-\hat
{d}_{0}\varepsilon)y_{1}^{e_{j}}+ \hat{d}_{0}\varepsilon,\varepsilon
y_{2}^{e_{j}} ) & \mathrm{in}\qquad(0,1)\times(-\frac{1}{2},\frac{1}{2})
,\vspace{0.1cm}\\
\mu_{\varepsilon}^{e_{j}}=0 & \mathrm{on}\qquad\partial((0,1)\times(-\frac
{1}{2},\frac{1}{2})) ,\vspace{0.1cm}\\
&
\end{array}
\right. \label{20}%
\end{equation}
Applying the result of \cite{Galdi}, Chap. III, p. 127 we get: there exist a
solution of \ref{20}) such that
\begin{equation}
\label{21}\|\mu_{\varepsilon}^{e_{j}} \|_{\left( H^{1}((0,1)\times(-\frac
{1}{2},\frac{1}{2}) )\right) ^{2}}=O(exp(-c/\varepsilon)).
\end{equation}
Expressing the norm of $\hat{U}^{e_{j}} $ with respect to the norm
$\mu_{\varepsilon}^{e_{j}} $ we obtain
\[
\|\hat{U}^{e_{j}} \|_{\left( H^{1}({\mathcal{B}}^{\varepsilon}_{j} \cap
\{x_{1}^{e_{j}}<\hat{d}_{0}\varepsilon\})\right) ^{2}}\leq\frac{1}%
{\varepsilon^{\frac{1}{2}}} \|\mu_{\varepsilon}^{e_{j}} \|_{\left(
H^{1}((0,1)\times(-\frac{1}{2},\frac{1}{2}) )\right) ^{2}} ,
\]
i.e. $\|\hat{U}^{e_{j}} \|_{\left( H^{1}({\mathcal{B}}^{\varepsilon}_{j}
\cap\{x_{1}^{e_{j}}<\hat{d}_{0}\varepsilon\})\right) ^{2}}%
=O(exp(-c/\varepsilon))$. So, $\Vert\hat{U}^{a} \Vert_{(H^{1}({\mathcal{B}%
}^{\varepsilon}))^{2}}=O(exp(-c/\varepsilon))$.

Define $U^{a}= \bar{u}^{a}-\hat{U}^{a}$. Then  $\left( U^{a},\bar{p}%
^{a}\right) $ satisfies the following problem~:
\begin{equation}
\left\{
\begin{array}
[c]{ll}%
-\mathrm{div}(\nu(x){\mathcal{D}}U^{a})+\nabla\bar{p}^{a}=f(x)+ \mathrm{div}%
(\nu(x){\mathcal{D}}\hat{U}^{a})+\Psi & \\
-\displaystyle{\sum_{e=e_{j};j=1,\dots,n}}\left\{  {\varepsilon}^{k+1}%
\Gamma_{j}\left(  F_{e_{j}}^{k}\right)  ^{T}\chi_{\varepsilon}(x)-\varepsilon
^{k}(\nabla_{\xi}\chi_{j}(\xi_{1}^{e_{j}})d_{k+1}^{e_{j}}-\nabla_{\tilde{\xi}%
}\chi_{j}(\tilde{\xi}_{1}^{e_{j}})\tilde{d}_{k+1}^{e_{j}})\right\}  , &
\text{in }{\mathcal{B}}^{\varepsilon},\\
\mathrm{div}\bar{u}^{a}=0 & \text{in }{\mathcal{B}}^{\varepsilon},\\
\bar{u}^{a}=g & \text{on }\partial{\mathcal{B}}^{\varepsilon},
\end{array}
\right.  \label{22}%
\end{equation}

\begin{theorem}
Let $(\bar{u}^{a},\bar{p}^{a})$ be the asymptotic solution given by
\eqref{15}, \eqref{16} and $(u_{\varepsilon},p_{\varepsilon})$ the solution of
\eqref{1}, the following estimates hold:%
\begin{equation}
\left\{
\begin{array}
[c]{l}%
\Vert\bar{u}^{a}-u_{\varepsilon}\Vert_{{H^{1}({\mathcal{B}}^{\varepsilon}%
)}^{2}}=O\left(  \varepsilon^{k+\frac{1}{2}}\right)  ,\\
\Vert\nabla\bar{p}^{a}-\nabla p_{\varepsilon}\Vert_{{H^{1}({\mathcal{B}%
}^{\varepsilon})}}=O\left(  \varepsilon^{k+\frac{1}{2}}\right)  .
\end{array}
\right.
\end{equation}

\end{theorem}

Applying a priori estimates and \eqref{21} we get
\[
\Vert U^{a}-u_{\varepsilon}\Vert_{H^{1}({\mathcal{B}}^{\varepsilon})}=O\left(
\varepsilon^{k+\frac{1}{2}}\right)
\]
and so,
\[
\Vert\bar{u}^{a}-u_{\varepsilon}\Vert_{H^{1}({\mathcal{B}}^{\varepsilon}%
)}=O\left(  \varepsilon^{k+\frac{1}{2}}\right)
\]
The estimate for the pressure is obtained from the a priori estimates. These
estimates justify the construction of the asymptotic expansion.

\begin{rque}
The main result can be easily generalized in case when the length of
$\beta_{j}$ is different from 1.
\end{rque}

\begin{rque}
Formula \eqref{11} shows that only $c^{e_{j}}_{0}$ could be different from
zero. The same analysis can be provided in the case of a multi-bundle
structure, that is the union of a finite number of thin domains of
${\mathcal{B}}^{\varepsilon}$ type (see \cite{pana2005}, section 4.5.2), and
in this case the constants $c^{e_{j}}_{l}$ should be determined from a system
of linear algebraic equations (see (4.5.43),(4.5.44) in \cite{pana2005}).
\end{rque}

\section{Numerical experiments}


\textbf{1.} We consider the Stokes flow in a rectangular domain $(0,1)\times
(0,\varepsilon)$ with $\varepsilon=0.1$:%

\begin{equation}
\label{N.1}\left\{
\begin{array}
[c]{ll}%
\displaystyle-\mathrm{div}(\nu(x_{1})Du)+\nabla p=0 & \displaystyle\text{ in
}(0,1)\times(0,\varepsilon)\\
\displaystyle\mathrm{div} u=0 & \displaystyle\text{ in }(0,1)\times
(0,\varepsilon)\\
\displaystyle u(0,x_{2})=\varepsilon^{2} \frac{x_{2}}{\varepsilon}\left(
1-\frac{x_{2}}{\varepsilon}\right)  & \displaystyle\text{ on }(0,\varepsilon
)\\
& \\
\displaystyle u(1,x_{2})=\varepsilon^{2} \frac{x_{2}}{\varepsilon}\left(
1-\frac{x_{2}}{\varepsilon}\right)  & \displaystyle\text{ on }(0,\varepsilon
)\\
\displaystyle u(x_{1},0)=u(x_{1},\varepsilon) & \displaystyle\text{ on }(0,1)
\end{array}
\right.
\end{equation}
\vspace{0.3cm} here $\nu(x_{1})=2x_{1}+2$. From \eqref{A.4} we get for $j=0$~:
\vspace{0.3cm}
\[
\left\{
\begin{array}
[c]{l}%
\displaystyle\frac{\partial}{\partial x_{1}}\left( \frac{1}{\nu(x_{1})}%
\frac{\partial q_{0}}{\partial x_{1}}\right) =0\\
\\
\displaystyle u_{0,1}=\varepsilon^{2}\frac{1}{\nu(x_{1})}\frac{\partial q_{0}%
}{\partial x_{1}}\frac{x_{2}}{\varepsilon}\left( 1-\frac{x_{2}}{\varepsilon
}\right) \\
\displaystyle u_{0,2}=0\\
\displaystyle p_{0}=0
\end{array}
\right.
\]
with the compatibility condition $\eqref{B.5}$ we have~:
\[
\int_{0}^{1}\xi_{2}(1-\xi_{2})-u_{0,1}(0,\xi_{2})\mathrm{d}\xi_{2}=0
\]

Since we assume that $p(1,x_{2})=0$ we get~: \vspace{0.3cm}
\[
\left\{
\begin{array}
[c]{l}%
\displaystyle q_{0}(x_{1})=-x_{1}(x_{1}+2)+3\\
\\
\displaystyle u_{0,1}=\varepsilon^{2}\frac{x_{2}}{\varepsilon}\left(
1-\frac{x_{2}}{\varepsilon}\right) \\
u_{0,2}=0\\
p_{0}=0
\end{array}
\right.
\]
\vspace{0.5cm} Solving numerically the problem \eqref{N.1} (by Comsol) we get
the following results~: for the first component of the velocity we have~:
\begin{figure}[h]
\begin{minipage}[h]{.48\linewidth}
\centering\epsfig{figure=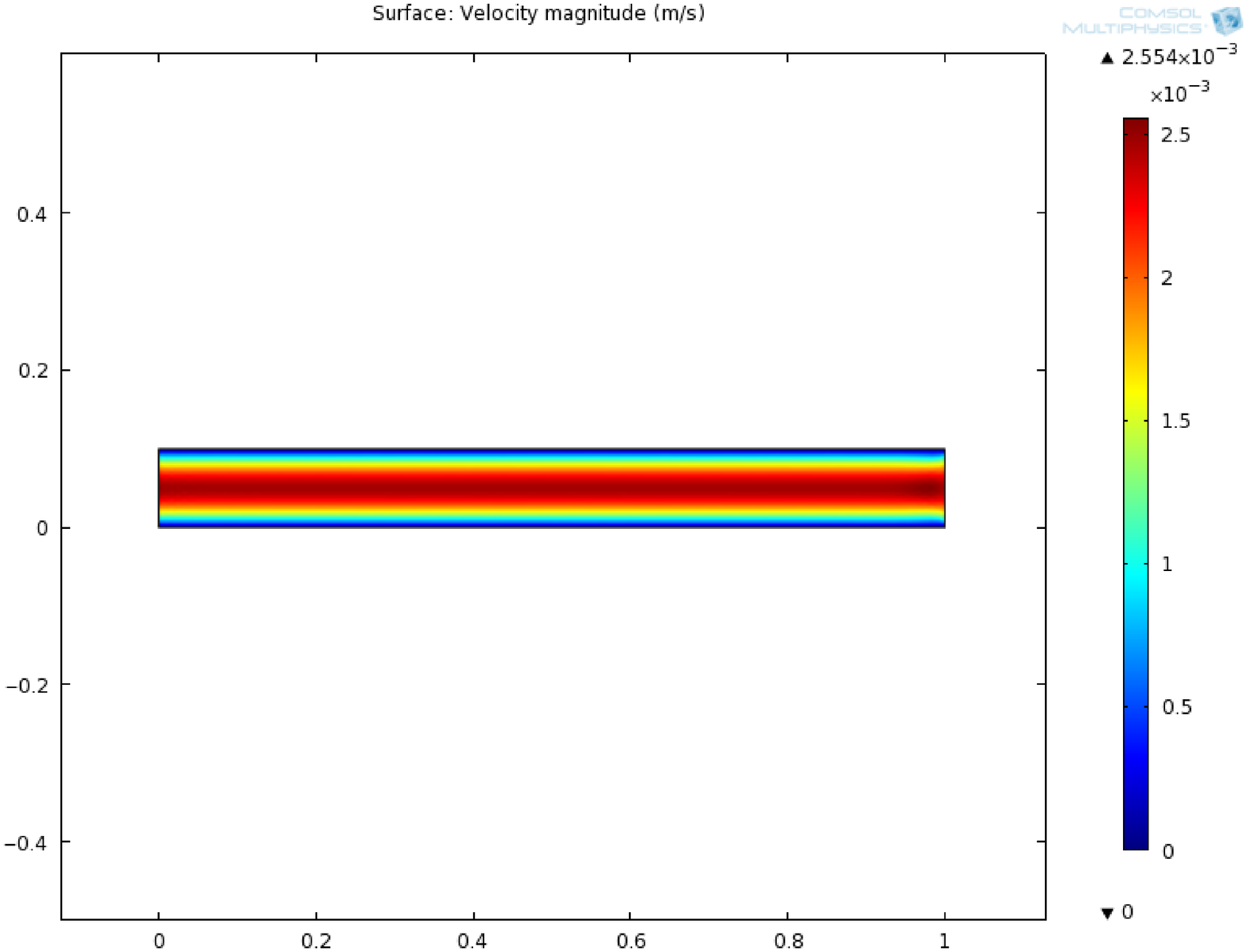,width=\linewidth}
\caption{First component of the velocity }\label{V1}
\end{minipage}\hfill\begin{minipage}[h]{.5\linewidth}
\centering\epsfig{figure=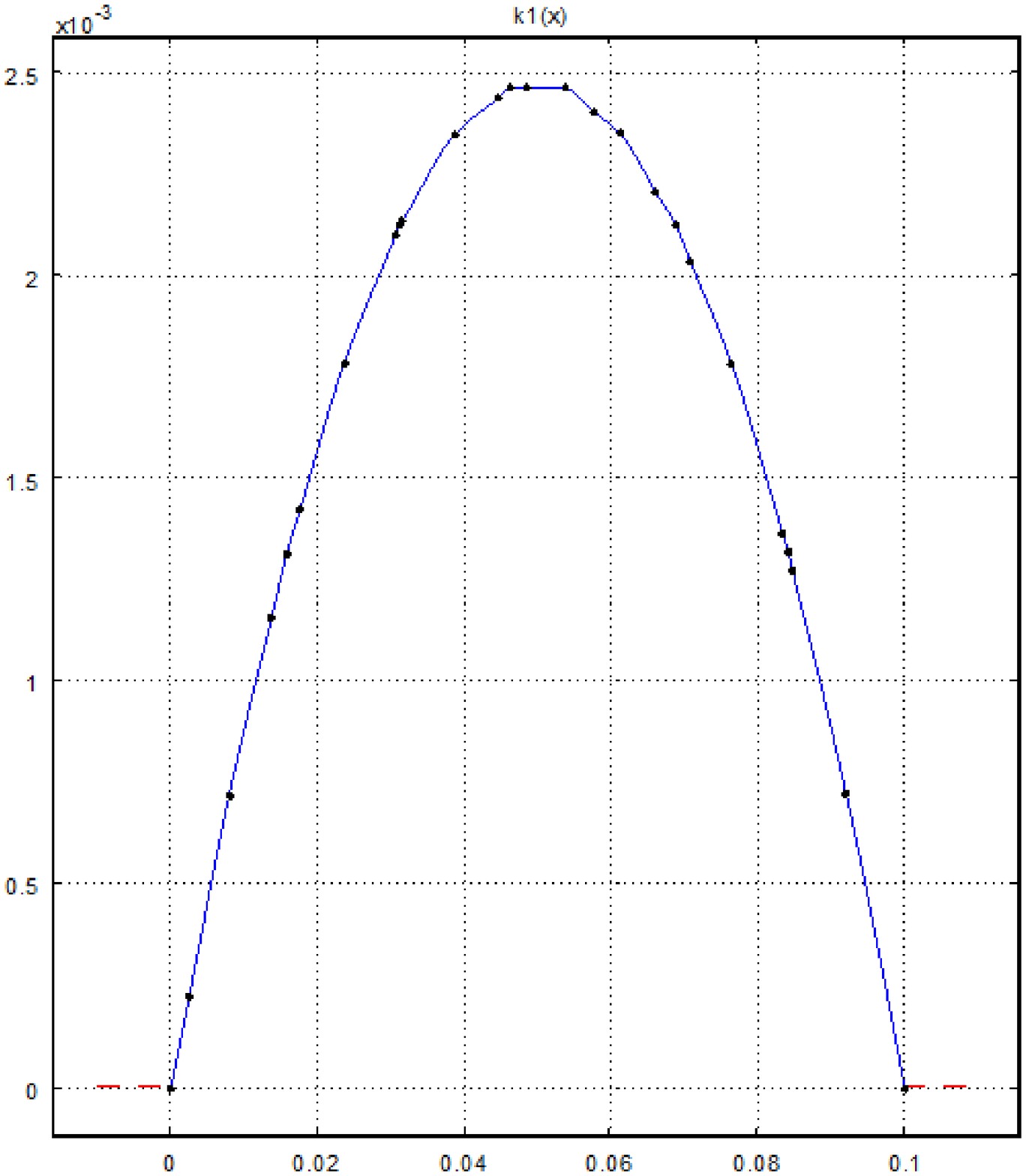,width=\linewidth}
\caption{The profile of the first component of the velocity for $x_1=0.5$ }\label{V2}
\end{minipage}
\end{figure}\clearpage
and for the pressure~:

\begin{figure}[h]
\begin{minipage}[h]{.48\linewidth}
\centering\epsfig{figure=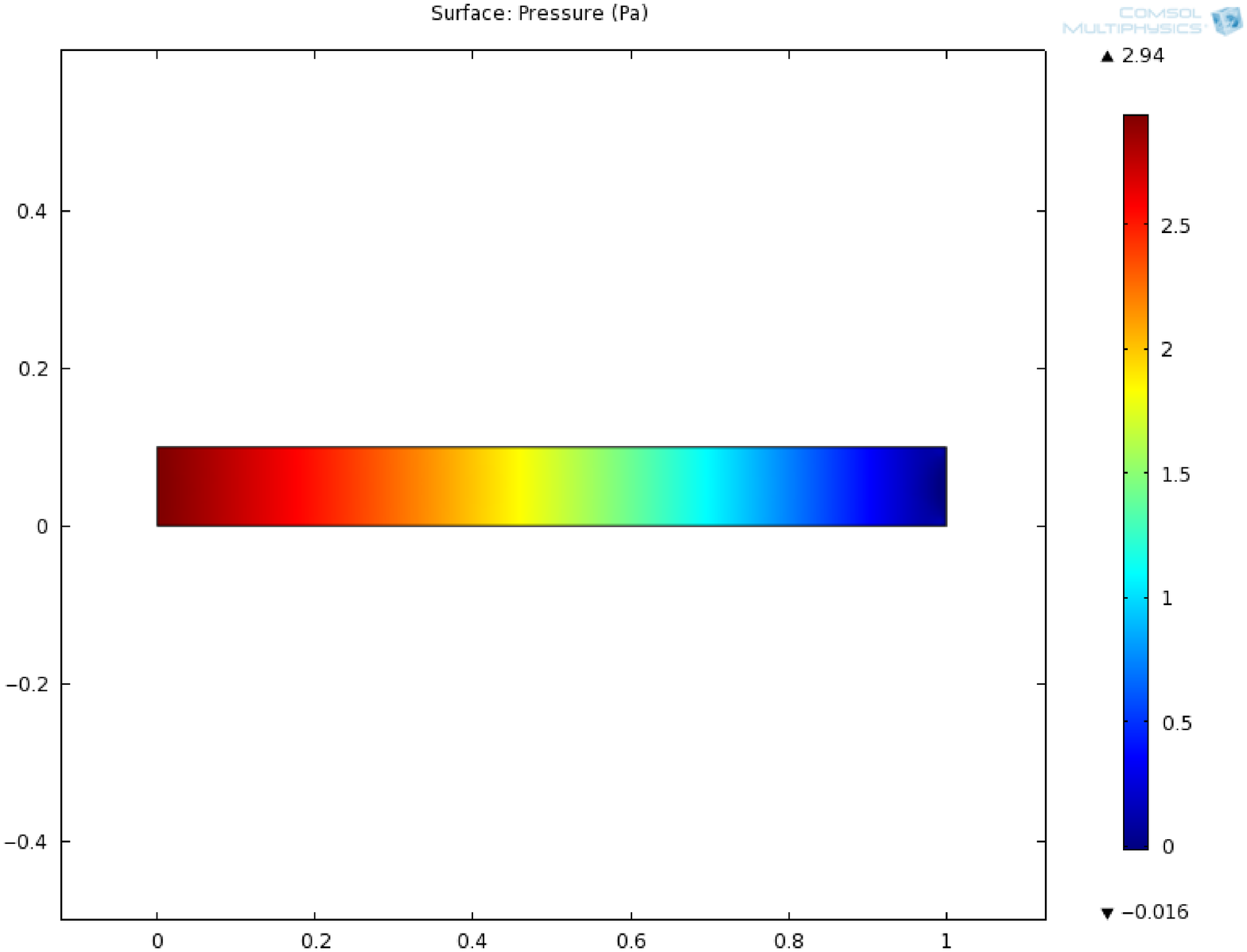,width=\linewidth}
\caption{The pressure }\label{P1}
\end{minipage}\hfill\begin{minipage}[h]{.5\linewidth}
\centering\epsfig{figure=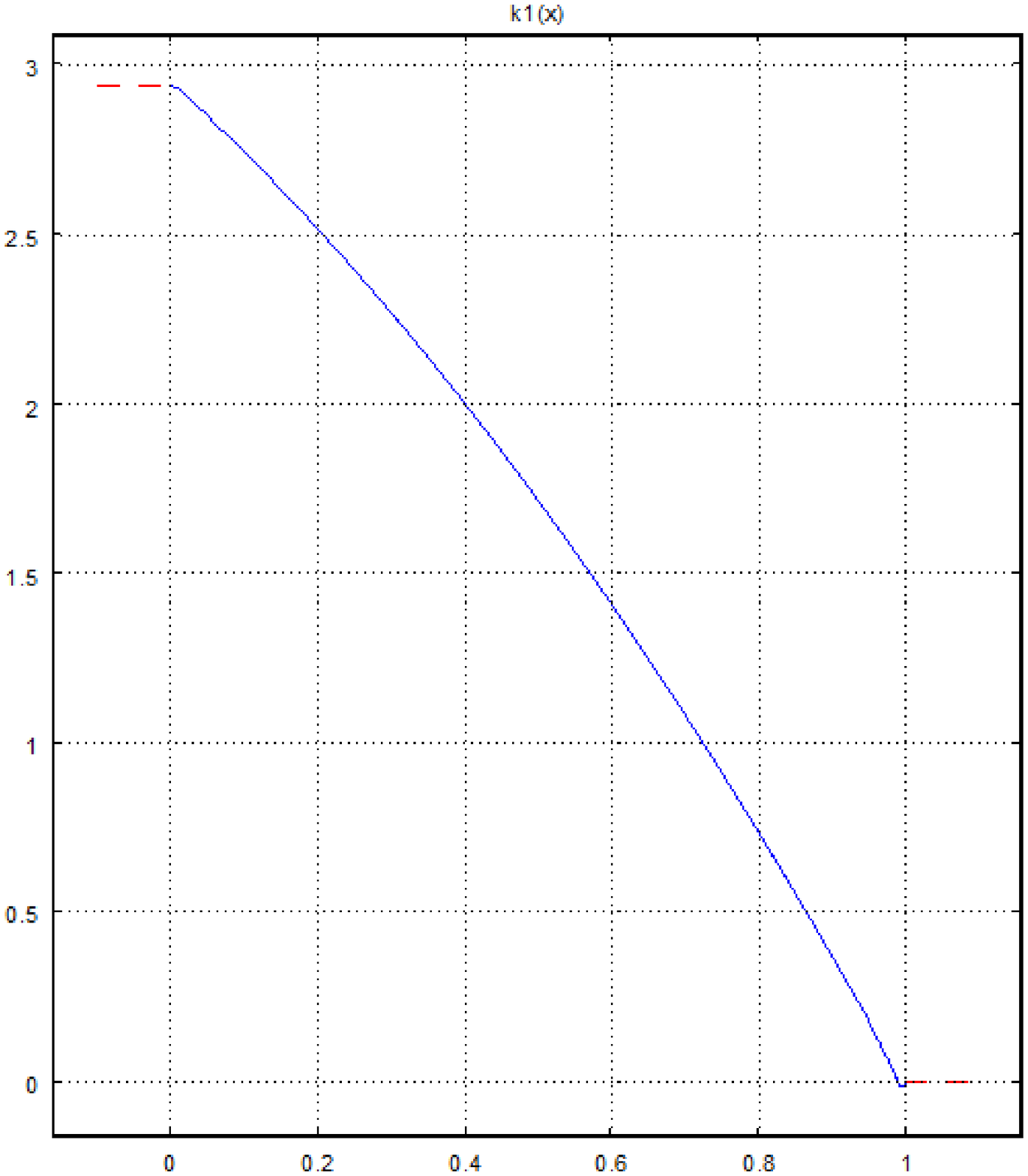,width=\linewidth}
\caption{Pressure profile in $x_2=0.05$ }\label{P2}
\end{minipage}
\end{figure}
The error between the numerical result and the leading term of the asymptotic
expansion is as follows~: \begin{figure}[h]
\begin{minipage}[h]{.48\linewidth}
\centering\epsfig{figure=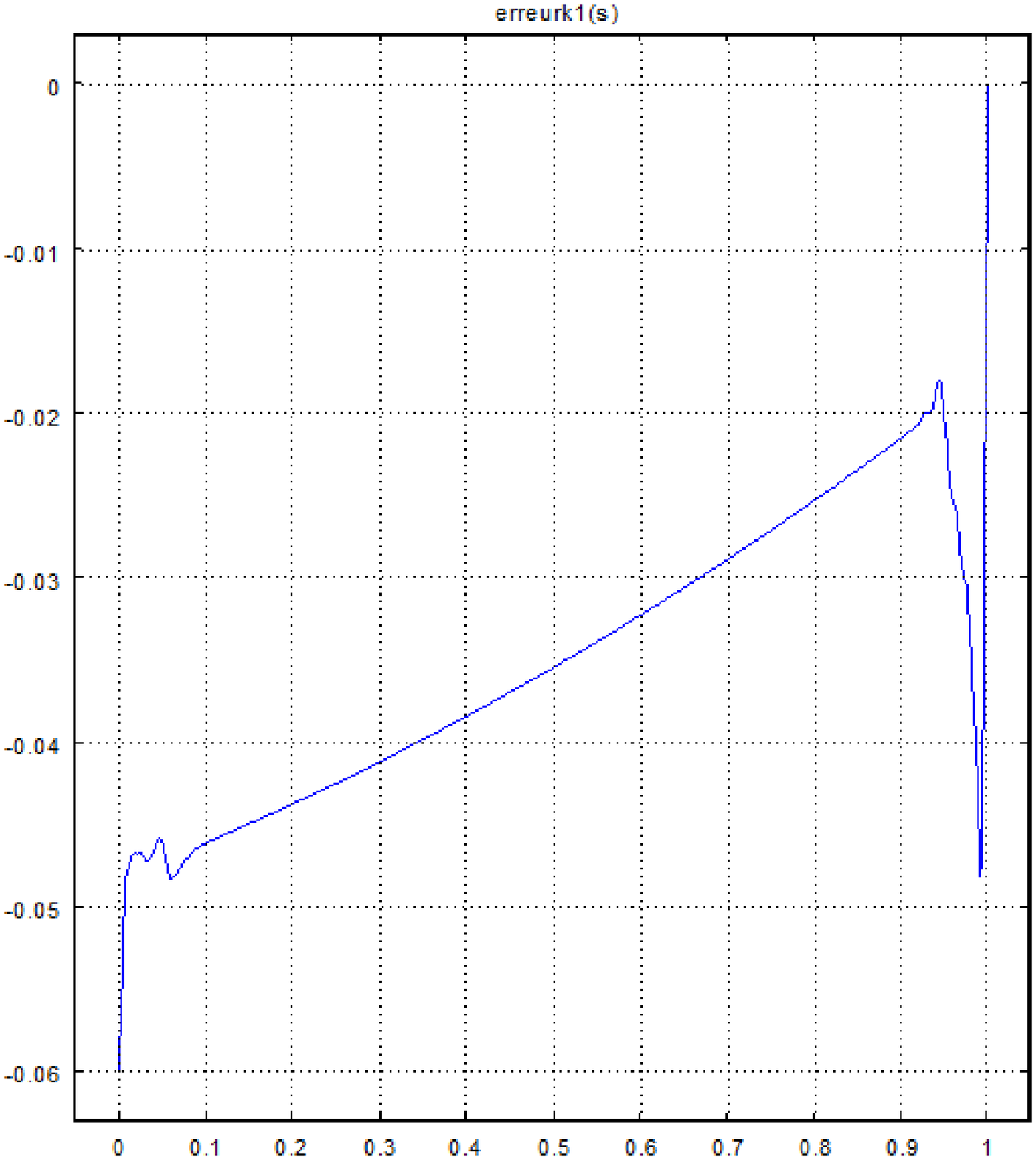,width=\linewidth}
\caption{Error for the pressure}\label{E1}
\end{minipage}\hfill\begin{minipage}[h]{.5\linewidth}
\centering\epsfig{figure=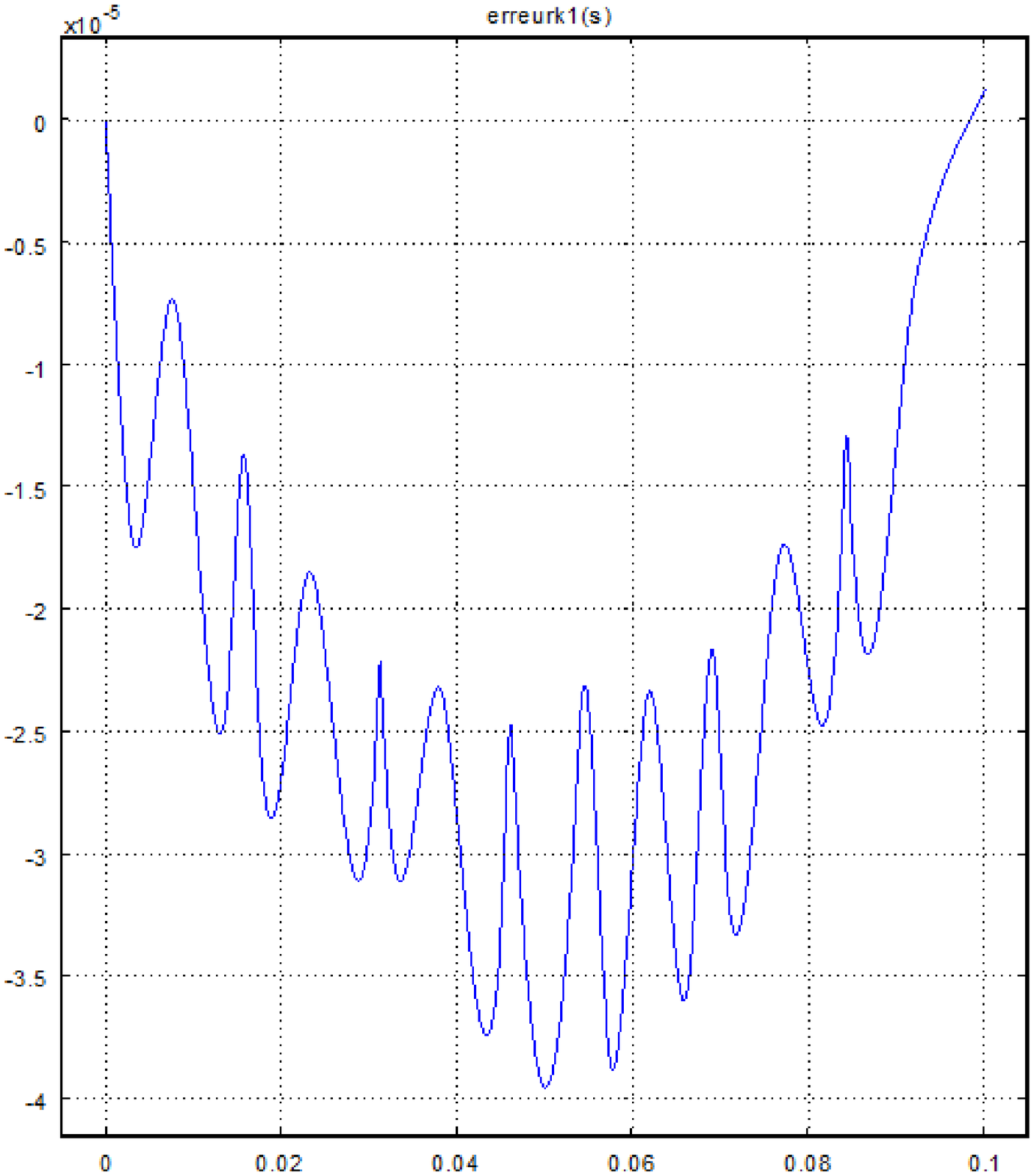,width=\linewidth}
\caption{Error for the velocity}\label{E2}
\end{minipage}
\end{figure}
We see that the error is of order of $\varepsilon$ et $\varepsilon^{5}$ for
the pressure and the velocity respectively.\newline\bigskip\textbf{2.}
Consider now the Stokes problem \eqref{N.1} in a T-shape domain ${\mathcal{B}%
}^{\varepsilon}=(-1,0]\times(0,\varepsilon)\cup(0,\varepsilon)\times
(-0.45,0.55)$
\begin{equation}
\left\{
\begin{array}
[c]{ll}%
\displaystyle u(1,x_{2})=g_{1}\left( \frac{x_{2}}{\varepsilon}\right)
=\varepsilon^{2} \frac{x_{2}}{\varepsilon}\left( 1-\frac{x_{2}}{\varepsilon
}\right)  & \text{ on }(0,\varepsilon)\\
& \\
\displaystyle u(x_{1},-0.45)=g_{2}\left( \frac{x_{1}}{\varepsilon}\right)
=\varepsilon^{2} \frac{x_{1}}{\varepsilon}\left( 1-\frac{x_{1}}{\varepsilon
}\right)  & \text{ on }(0,\varepsilon)\\
& \\
\displaystyle u(x_{1},-0.45)=g_{3}\left( \frac{x_{1}}{\varepsilon}\right)
=2\varepsilon^{2} \frac{x_{1}}{\varepsilon}\left( 1-\frac{x_{1}}{\varepsilon
}\right)  & \text{ on }(0,\varepsilon)\\
& \\
&
\end{array}
\right.
\end{equation}
as inflow/outflow conditions, $\varepsilon=0.1$, $u=0$ anywhere else on the
boundary and $\displaystyle\nu(x_{1},x_{2})=2+2(x_{1}+1)\mu_{1}(x_{1}%
)+2(x_{2}+1)\mu_{2}(x_{2})+2(x_{2}+1)\mu_{3}(x_{2})$. $\mu_{i}$ are defined in
such a way that $\nu$ is equal to $2x_{1}^{e_{j}}+2$ in each rectangle and
constant near the nodes. \clearpage
The functions $\mu_{1},\mu_{2}$ and $\mu_{3}$ have the following shape~:
\begin{figure}[th]
\begin{minipage}[h]{.32\linewidth}
\centering\epsfig{figure=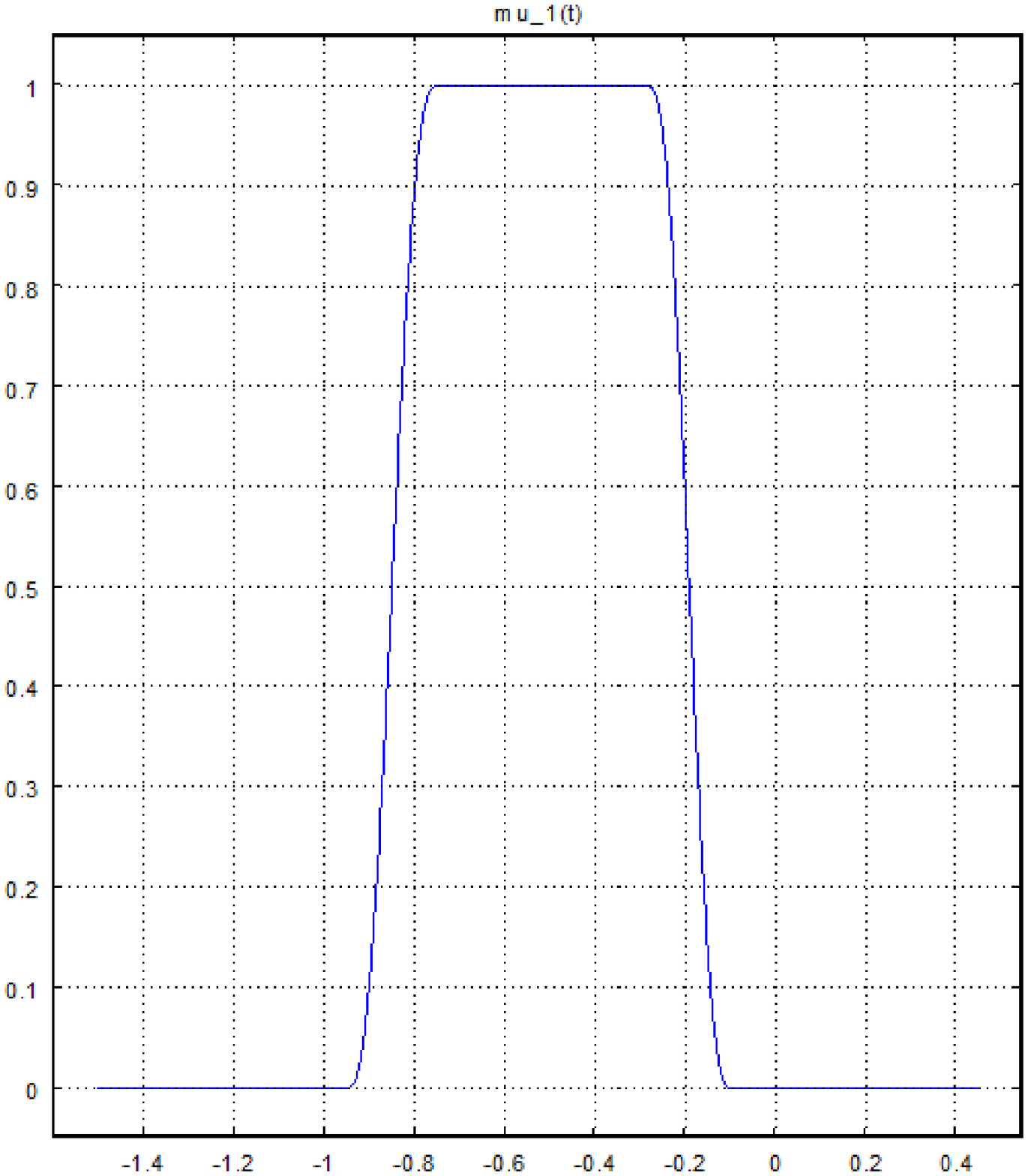,width=\linewidth}%
\caption{$\mu_1$}
\end{minipage}\hfill\begin{minipage}[h]{.32\linewidth}
\centering\epsfig{figure=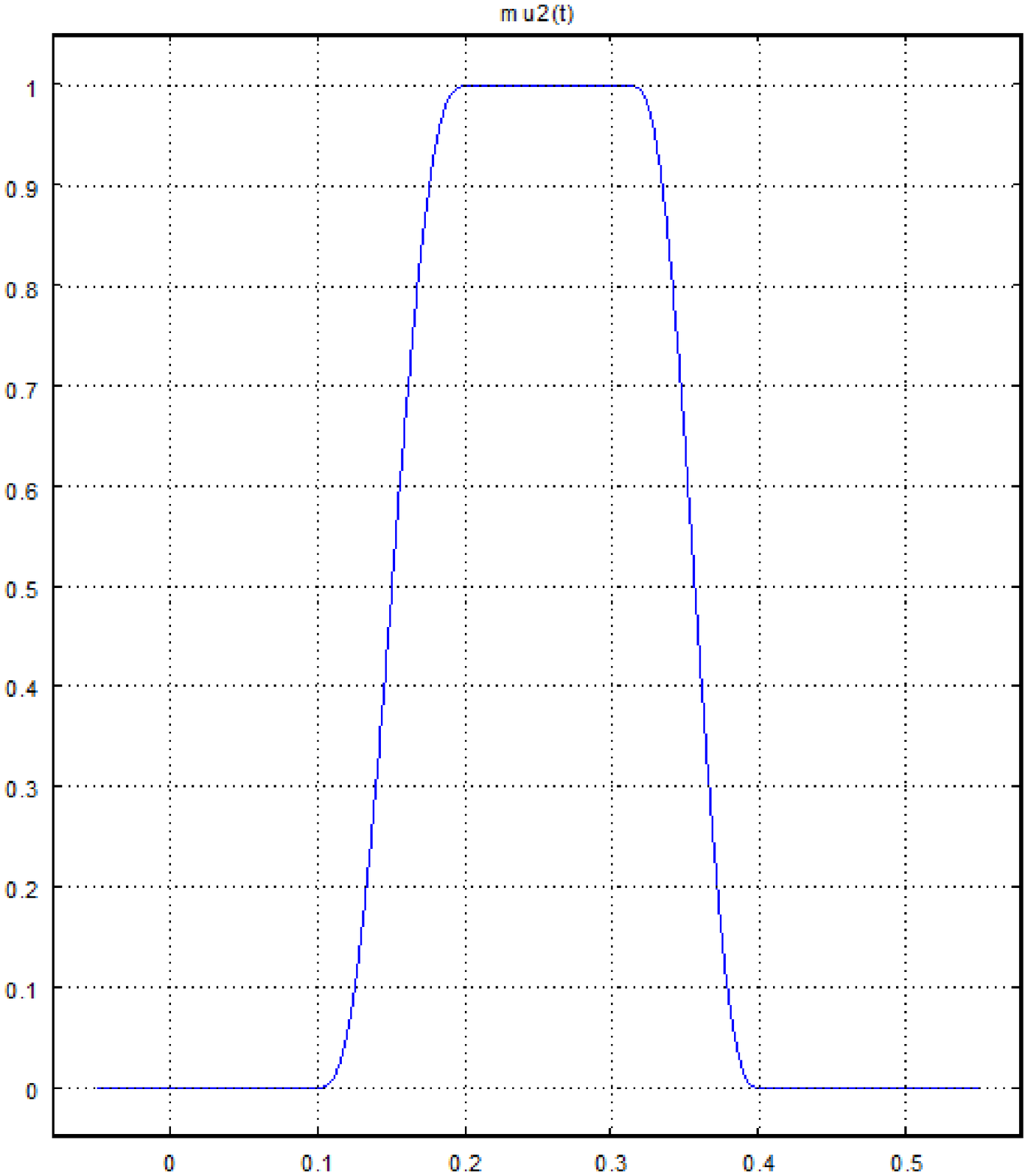,width=\linewidth}%
\caption{$\mu_2$}
\end{minipage}\hfill\begin{minipage}[h]{.32\linewidth}
\centering\epsfig{figure=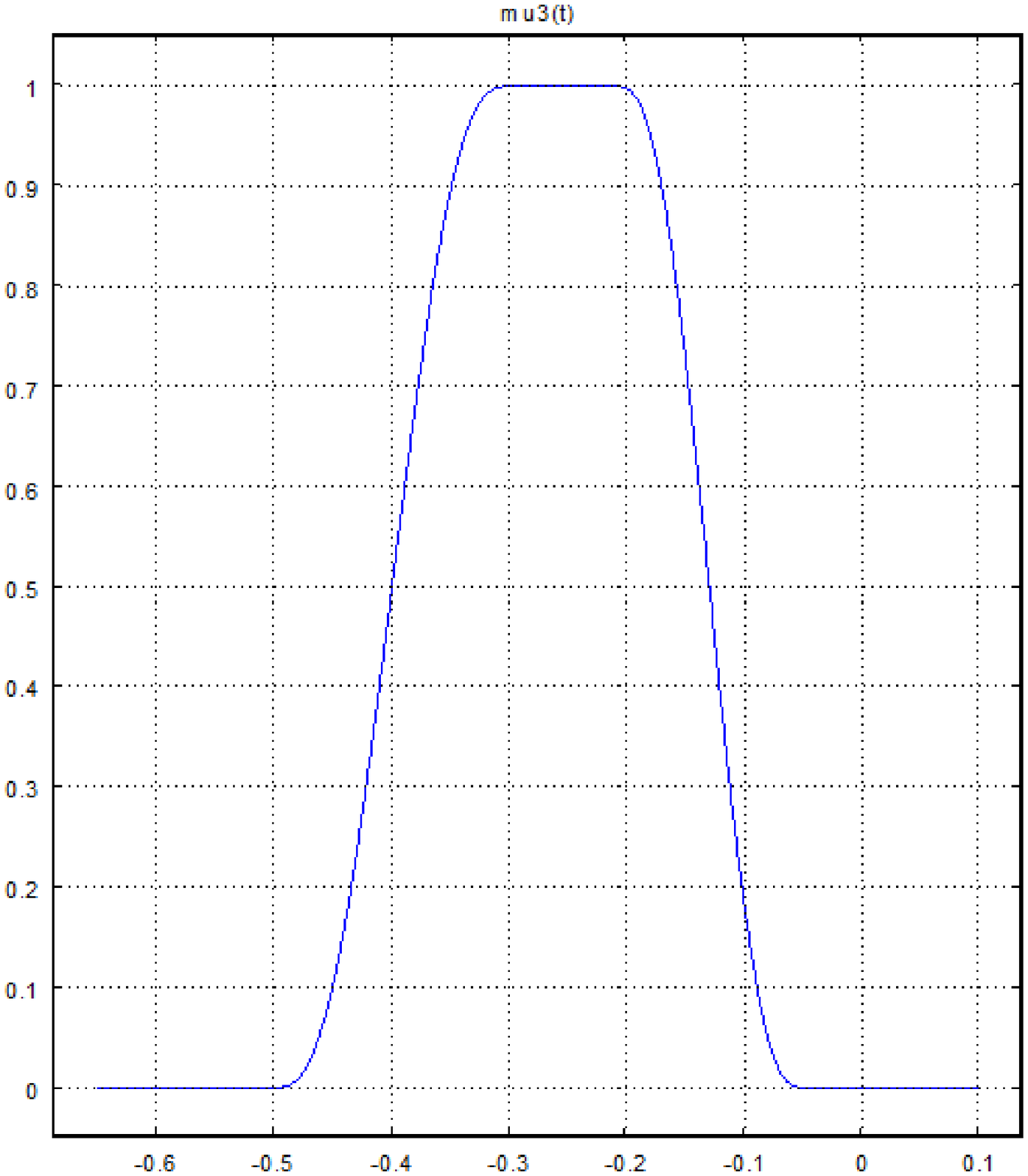,width=\linewidth}%
\caption{$\mu_3$}
\end{minipage}
\end{figure}\newline So we get \begin{figure}[th]
\begin{minipage}[h]{.48\linewidth}
\centering\epsfig{figure=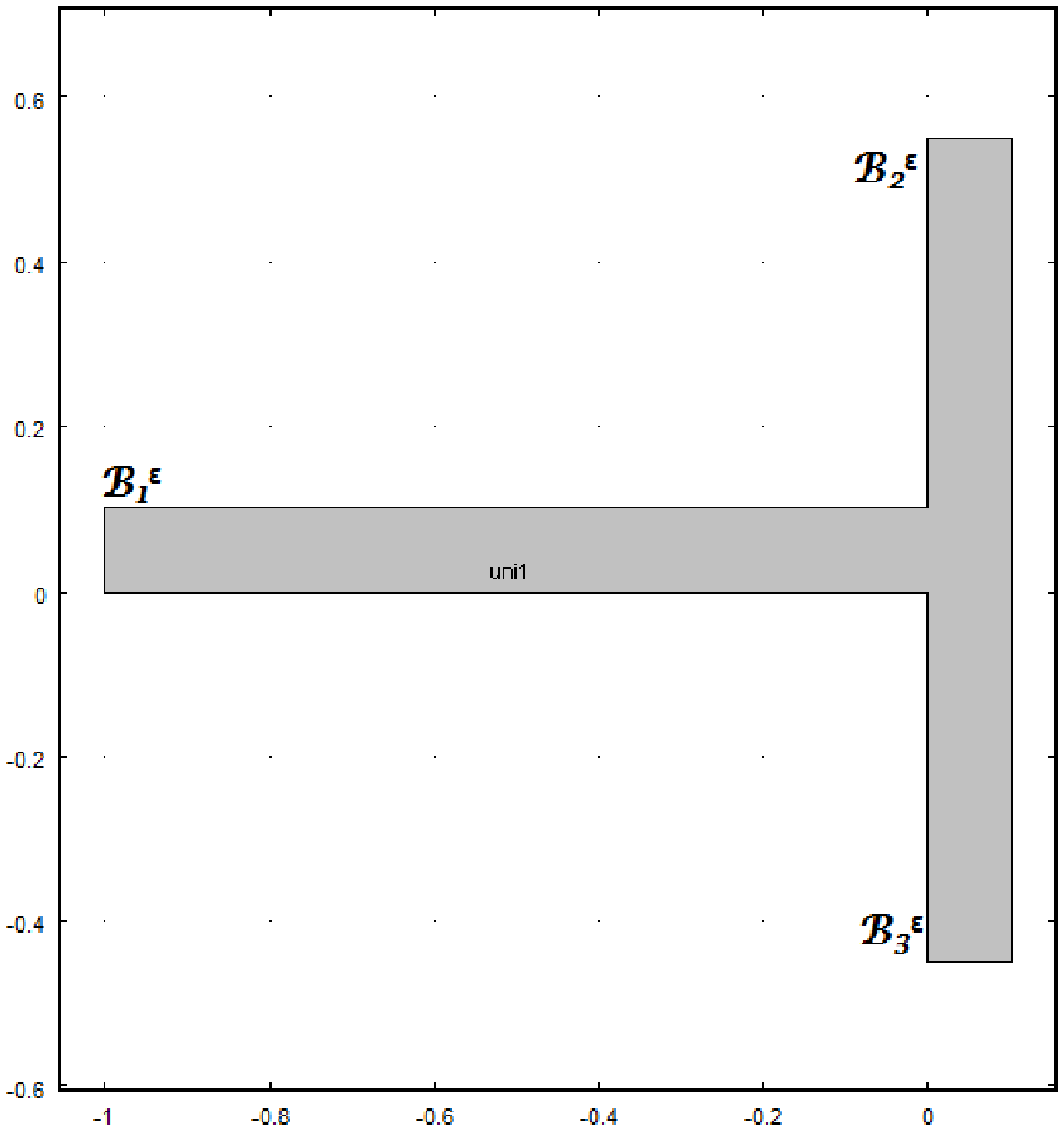,width=\linewidth}%
\caption{Thin domain $\B^\varepsilon$}
\end{minipage}\hfill\begin{minipage}[h]{.48\linewidth}
\centering\epsfig{figure=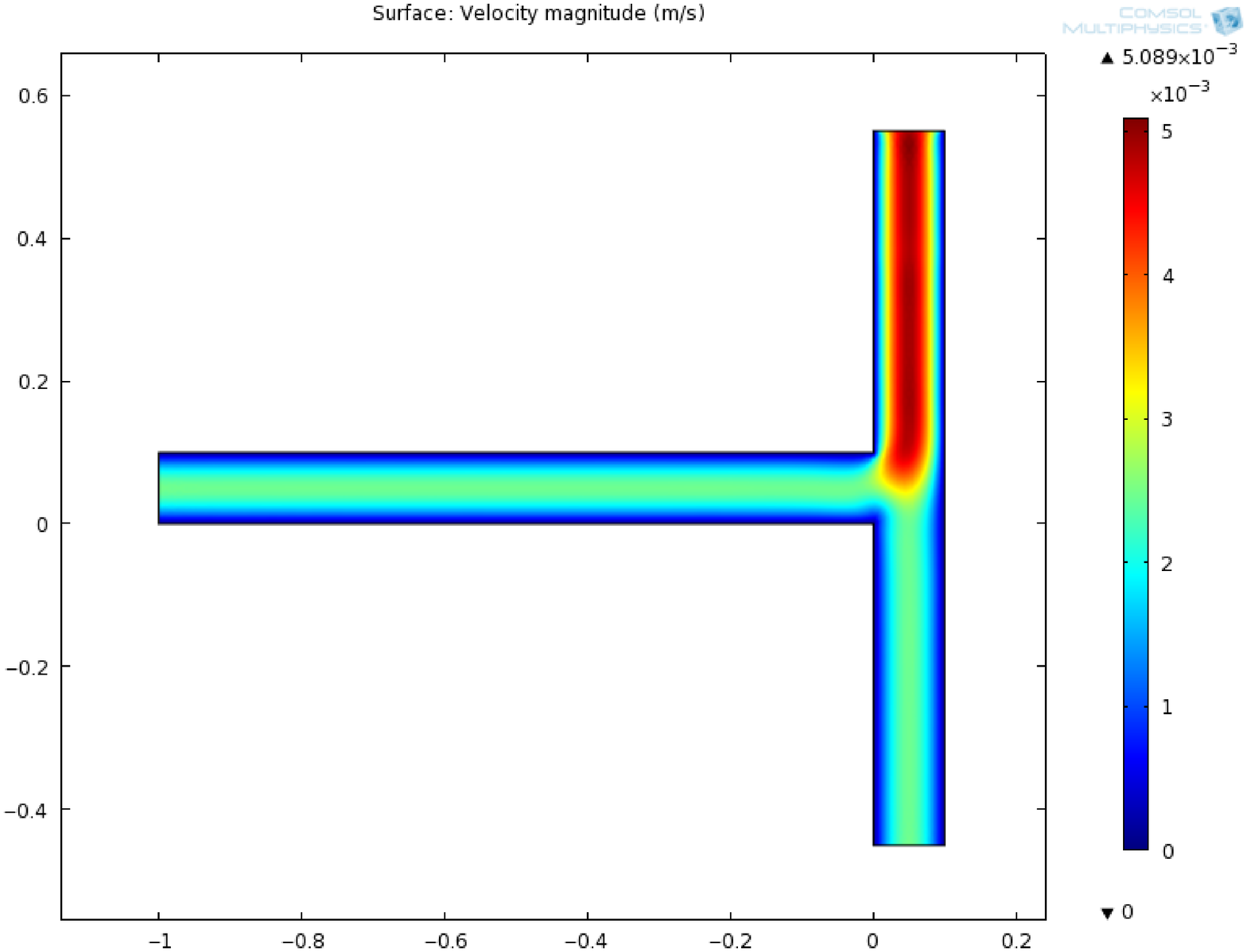,width=\linewidth}%
\caption{Velocity magnitude}
\end{minipage}
\end{figure}\newline\vspace{0.2cm} We do two crop sections, the first in
${\mathcal{B}}_{2}^{\varepsilon}=(0,\varepsilon)\times(0,0.55)$ and the second
in ${\mathcal{B}}_{3}^{\varepsilon}=(0,\varepsilon)\times(-0.45,0)$, and the
results are as follows~: \vspace{0.2cm} \begin{figure}[th]
\begin{minipage}[h]{.48\linewidth}
\centering\epsfig{figure=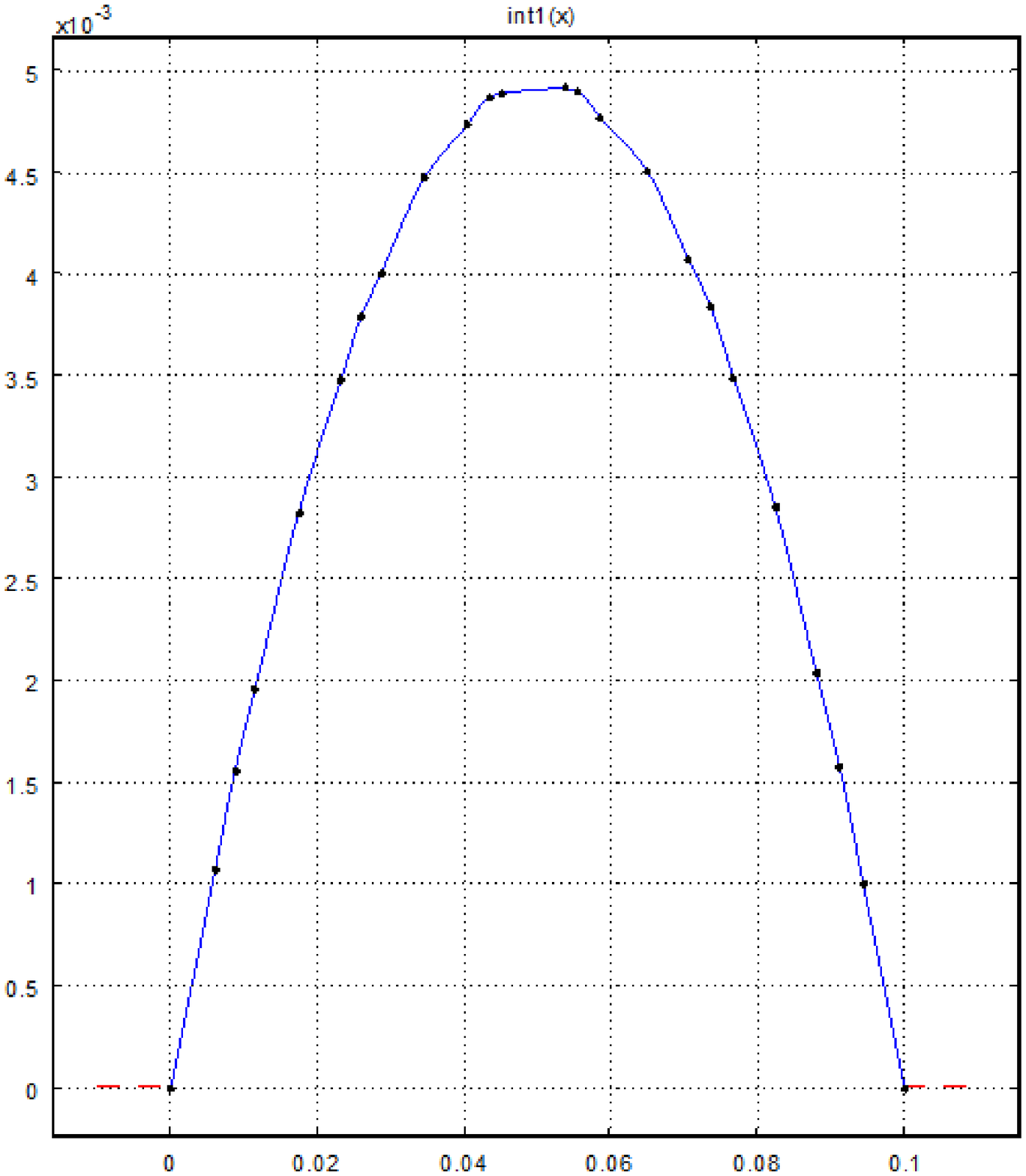,width=\linewidth}
\caption{Profile of the first component of the velocity in $\B_2^\varepsilon$ for $x_2=0.2$}\label{MT2}
\end{minipage}\hfill\begin{minipage}[h]{.5\linewidth}
\centering\epsfig{figure=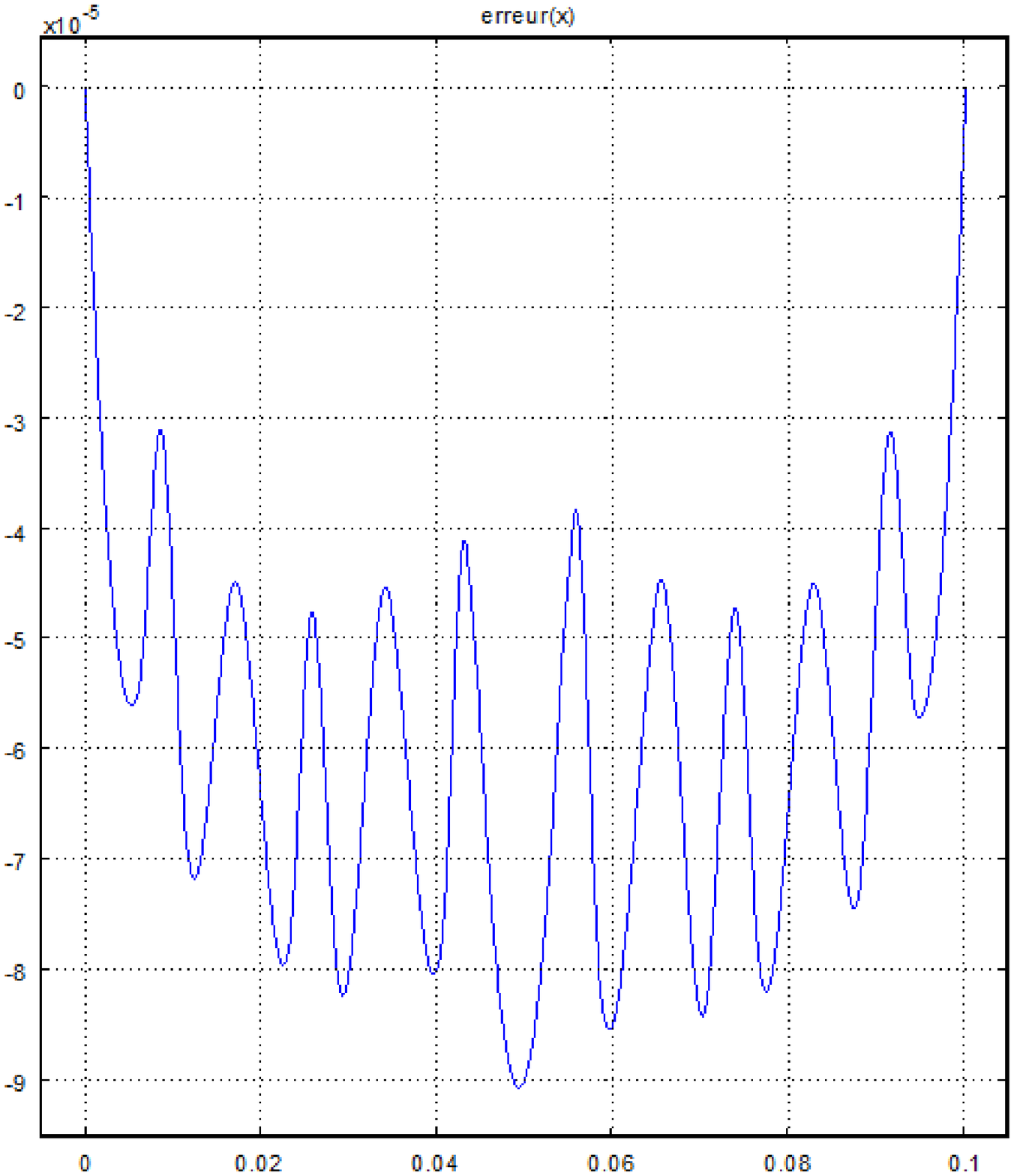,width=\linewidth}
\caption{Error for the velocity in $\B_2^\varepsilon$}\label{MT3}
\end{minipage}
\end{figure}\begin{figure}[h]
\begin{minipage}[h]{.48\linewidth}
\centering\epsfig{figure=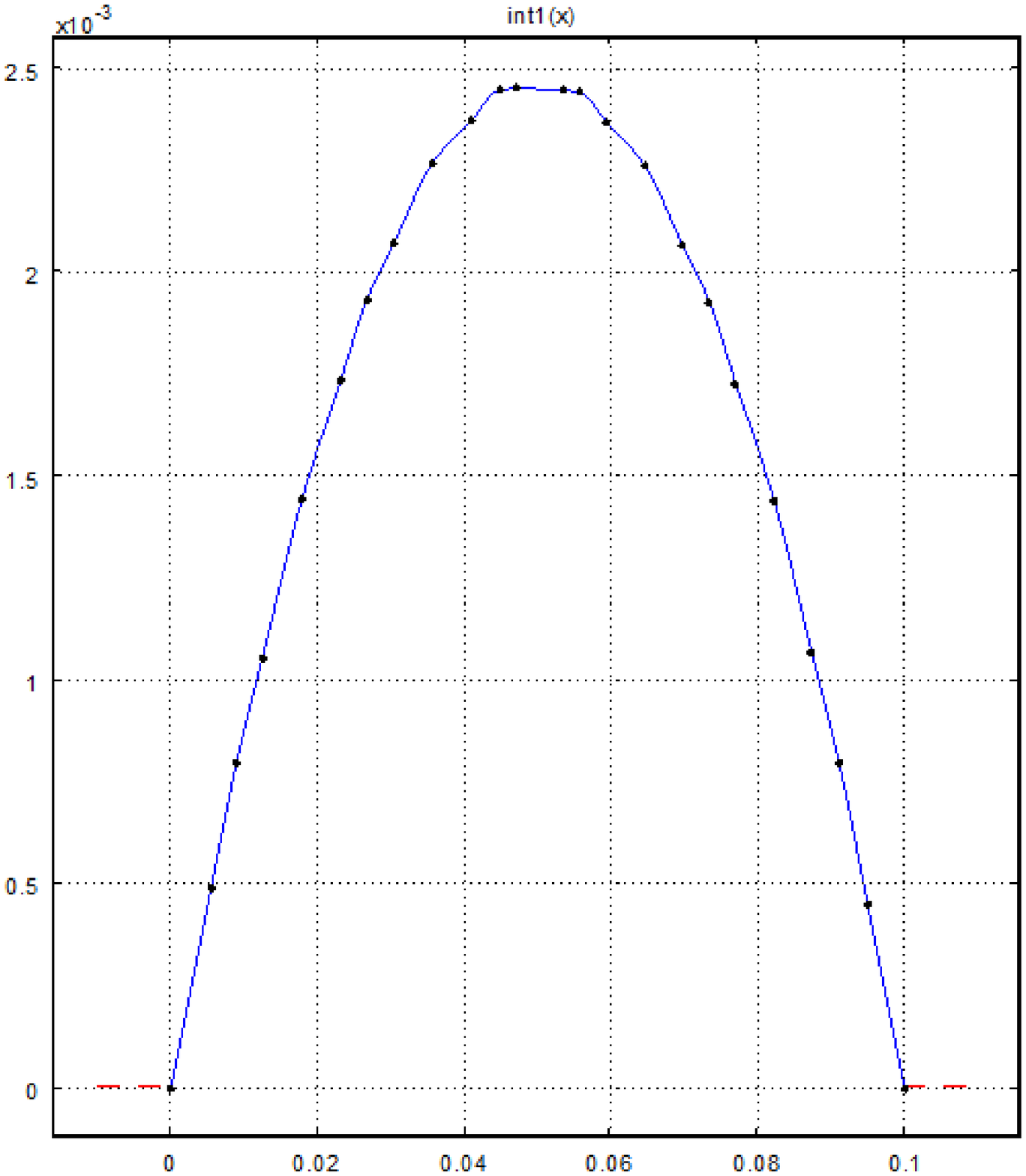,width=\linewidth}
\caption{Profile of the first component of the velocity in $\B_3^\varepsilon$ for $x_2=-0.2$}\label{MT4}
\end{minipage}\hfill\begin{minipage}[h]{.5\linewidth}
\centering\epsfig{figure=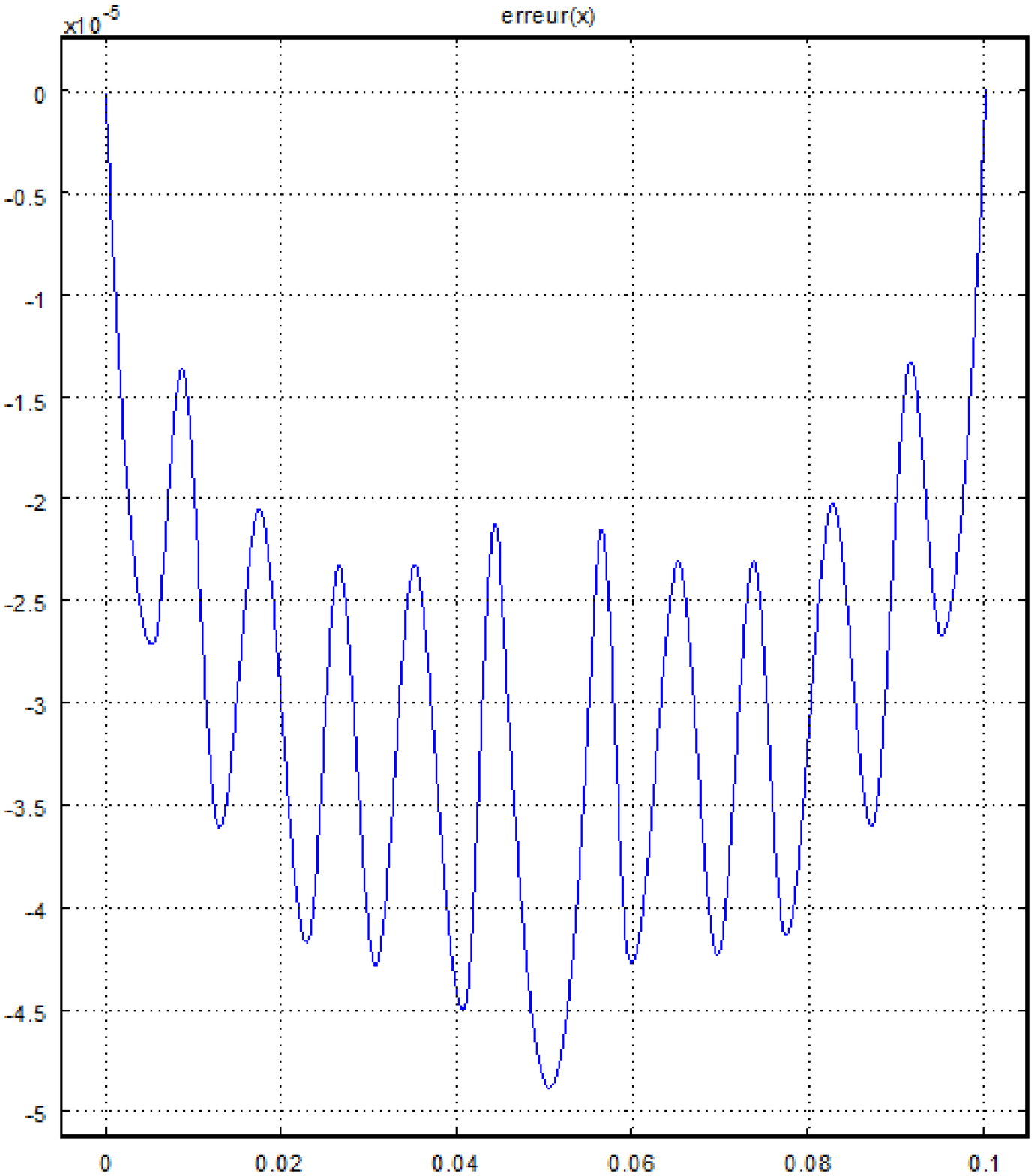,width=\linewidth}
\caption{Error for the velocity in $\B_3^\varepsilon$}\label{MT5}
\end{minipage}
\end{figure}\newpage We see that the error is of order of $\varepsilon^{5}$ as
in the case of a rectangle and estimates. It confirms the theoretical
prediction in section3.\newline

\bigskip\textbf{Acknowledgments.}{ The authors were partially supported by the
following grants: "Strutture sottili" of the program "Collaborazioni
interuniversitarie internazionali" (2004-2006) of the Italian Ministry of
Education, University and Research; SFR MOMAD of the university of Saint
Etienne and ENISE ( the Ministry of the Research and Education of France), the
joint French-Russian PICS CNRS grant "Mathematical modeling of blood
diseases"  and by the grant no. 14.740.11.0875 "Multiscale problems: analysis
and methods" of the Ministry of Edication and Research  of Russian Federation.
}

\end{document}